%%%%%%%%%%%%%%%%%%%%%%%%%%%%%%%%

\documentstyle{amsppt}

\loadbold

\NoRunningHeads \magnification=\magstephalf
\hcorrection{.5truein}
\NoBlackBoxes \refstyle{C} \loadbold

  \overfullrule=0pt

\def\thfill{\null\nobreak\hfill}
\def\endbeweis{\thfill\vbox{\hrule
   \hbox{\vrule\hbox to 5pt{\vbox to 5pt{\vfil}\hfil}\vrule}\hrule}}

\def\rtimes{\times\hskip-4.1pt\hbox{\vrule height 4.59 pt width 0.04
em}\ }

\def\ref#1{\par\noindent\hangindent3\parindent
\hbox to 3\parindent{#1\hfil}\ignorespaces}

\topmatter

\title  On the existence of $F$-crystals\endtitle

   \author
R. Kottwitz \footnote{ Mathematics Department, University of
Chicago, 5734
University Ave, Chicago, IL  60637, USA}
   and M. Rapoport
\footnote{ Mathematisches Institut der Universit\"at zu K\"oln,
Weyertal
86-90, D--50931 K\"oln, Germany}
   \endauthor

\endtopmatter

\document

\bigskip
By an $F$-isocrystal we mean a pair $(N,F)$, consisting of a
finite-dimensional vector space over the fraction field $L$ of the ring
$W(\overline{\bold F}_p)$ of Witt vectors of $\overline{\bold F}_p$
and a
Frobenius-linear bijective endomorphism of $N$. Isocrystals form a category
in an
obvious way. By Dieudonn\'e, $F$-isocrystals are classified up to
isomorphism by their Newton slope sequence. More precisely, let
$$({\bold
Q}^n)_+= \{  (\nu_1,\ldots,\nu_n)\in{\bold Q}^n;\ \nu_1\geq
\nu_2\geq\ldots\geq \nu_n\}\ \ .$$ Then we obtain an injective map
(the
Newton map) $$\{ \hbox{isocrystals of dimension}\
n\}/\simeq\longrightarrow ({\bold Q}^n)_+\ ,\ (N,F)\longmapsto
\boldsymbol\nu (N,F)\ \ .$$ Its image is characterized by the
following integrality condition. Let us write $\nu\in ({\bold Q}^n)_+$
in
the form $$\nu=(\nu(1)^{m_1},\ldots , \nu(r)^{m_r})\ \ ,\ \
\hbox{where}\
\nu(1)>\nu(2) >\ldots >\nu(r)\ \ .$$ Then the integrality condition
states
that $m_i\nu(i)\in {\bold Z}$, $\forall i=1,\ldots r$.
\par
Let now $(N,F)$ be an isocrystal of dimension $n$. Let $M$ be a
$W(\overline{\bold F}_p)$-lattice in $N$. Then the relative position of
$M$ and $FM$ is measured by the Hodge slope sequence
$\boldsymbol\mu=\boldsymbol\mu(M)={\roman{inv}}(M,FM)\in
({\bold Z}^n)_+$. Here $({\bold Z}^n)_+={\bold Z}^n\cap ({\bold
Q}^n)_+$,
and $(\mu_1,\ldots, \mu_n)\in ({\bold Z}^n)_+$ equals
$\boldsymbol\mu(M)$ iff there exists a $W(\overline{\bold
F}_p)$-basis $e_1,\ldots, e_n$ of $M$ such that
$p^{\mu_1}e_1,\ldots,
p^{\mu_n}e_n$ is a $W(\overline{\bold F}_p)$-basis of $F(M)$.
Mazur's
inequality states that $$\boldsymbol\mu (M)\geq
\boldsymbol\nu(N,F)\ \ ,$$ where the partial order relation on
$({\bold Q}^n)_+$ is the usual dominance order, comp.\ section 1.
\par
One result in this paper is a converse to this statement.
\medskip\noindent
{\bf Theorem A:} {\it Let $(N,F)$ be an isocrystal of dimension $n$.
Let
$\mu\in ({\bold Z}^n)_+$ be such that $\mu\geq
\boldsymbol\nu(N,F)$.
Then there exists a $W(\overline{\bold F}_p)$-lattice $M$ in $N$
with
$\mu=\boldsymbol\mu(M)$.}
\medskip
This is the content of Theorem 4.11 below which also gives the
corresponding statement for the group of symplectic similitudes.
As a
matter of fact, one can formulate a corresponding statement for any
quasi-split group over a $p$-adic field $F$ which splits over an
unramified extension of $F$ ([R], section 4), and we
conjecture that
Theorem A is true in this generality. Much of our argument in
section 4
below is formulated in the context of a split group with simply
connected
derived group, but we have not carried out the proof in this
generality.
Also note that if $\boldsymbol\nu(N,F)\in ({\bold Z}^n)_+$, the
general case was handled in [R] as an application of the
positivity
property of the Satake isomorphism. This positivity property also
plays a
crucial role in our proof of Theorem A. We also note that when
$\boldsymbol\nu (N,F)$ is of the form $\boldsymbol\nu
(N,F)=(\nu,\ldots, \nu)$, one can write down explicitly a lattice $M$
as
in Theorem A, and similarly in the more general case when $\mu$
is
decomposable with respect to $\boldsymbol\nu (N,F)$ (i.e., the
Hodge
polygon passes through all break points of the Newton polygon).
The
general case is reduced to this decomposable case, but then it
does not
seem so easy to produce explicitly a lattice $M$ with the required
properties.
\par
Theorem A may be considered as a statement on generalized affine
Deligne-Lusztig varieties. Let $\bar I\subset {\bold Z}/n{\bold Z}$ be
a
non-empty subset and let $M_{\bullet}$ be a periodic lattice chain
of type
$\bar I$. Then the relative position of $M_{\bullet}$ and
$FM_{\bullet}$
is an element $\boldsymbol\mu(M_{\bullet})\in \tilde W^{\bar
I}\setminus \tilde W / \tilde W^{\bar I}$. Here $\tilde W={\bold
Z}^n\rtimes S_n$ is the extended affine Weyl group of
${\roman{GL}}_n$ and
$\tilde W^{\bar I}$ is the parabolic subgroup of $\tilde W$
corresponding
to $\bar I$. The {\it generalized affine Deligne-Lusztig variety of type
$\bar I$ corresponding to} $w\in \tilde W^{\bar I}\setminus \tilde W
/\tilde W^{\bar I}$ is the set of all periodic lattice chains
$M_{\bullet}$ of type $\bar I$ with $\boldsymbol\mu(M_{\bullet})=
w$
(comp.\ [R], section 4). It seems a difficult question to
determine
for which $w$ this set is non-empty. Theorem A gives an answer to
this
question in case $\bar I= \{ 0\}$, in which case a periodic lattice
chain
of type $\bar I$ is simply a lattice and $\tilde W^{\bar I}\setminus
\tilde W /\tilde W^{\bar I}$ can be identified with $({\bold Z}^n)_+$.
\par
The question raised above becomes more tractable in case we
form a certain
finite union of Deligne-Lusztig varieties. Let $\mu\in ({\bold Z}^n)_+$
be
a minuscule element (i.e.\ $\mu_1-\mu_n\leq 1$) and consider
${\bold Z}^n$
as a subgroup of $\tilde W$. Let $${\roman{Adm}}(\mu)=\{ w\in
\tilde W;\
w\leq \mu'\ \ \hbox{for some}\ \mu'\in S_n\mu\}$$ be the {\it
$\mu$-admissible set} ([KR]). For a non-empty subset $\bar I$ let
${\roman{ Adm}}_{\bar I}(\mu)$ be the image of
${\roman{Adm}}(\mu)$ in
$\tilde W^{\bar I}\setminus \tilde W /\tilde W^{\bar I}$. We note that
by
[KR] this coincides with the {\it $\mu$-permissible subset} of $\tilde
W^{\bar I}\setminus \tilde W/\tilde W^{\bar I}$. Let $X(\mu, F)_{\bar
I}$
be the union of the generalized Deligne-Lusztig varieties of type
$\bar I$
corresponding to elements in ${\roman{Adm}}_{\bar I}(\mu)$.
\par
Our second main result in this paper is the following theorem.
\medskip\noindent
{\bf Theorem B:} {\it Let $\mu= \omega_r=(1^r, 0^{n-r})$ for some
$0\leq
r\leq n$.
\medskip
\item{(i)} For any non-empty subset $\bar I\subset {\bold Z}/n{\bold
Z}$,
$$X(\mu, F)_{\bar I}\neq\emptyset\ \hbox{if and only if}\ \mu\geq
\boldsymbol\nu(N,F)$$
\item{(ii)} For any non-empty subsets $\bar I$ and $\bar J$ of
${\bold Z}/n{\bold Z}$
with $\bar J\subset \bar I$ the forgetful map $$X(\mu,F)_{\bar
I}\longrightarrow X(\mu, F)_{\bar J}$$
\item{}is surjective.}
\medskip
This is the content of Proposition 1.1, which concerns the group
${\roman{GL}}_n$. Proposition 2.1 is the analogous statement for
the group
${\roman{GSp}}_{2n}$, i.e.\ for isocrystals with a symplectic
structure.
In section 3 we formulate the general problem. Section 4 is devoted
to the
proof of Theorem A, for ${\roman{GL}}_n$ and
${\roman{GSp}}_{2n}$. In
section 5 we treat the groups $R_{F'/F}{\roman{GL}}_n$ and
$R_{F'/F}{\roman{GSp}}_{2n}$ (restriction of scalars from a finite
unramified extension). If we had proved Theorem A for all unramified
reductive groups, this section could be eliminated. In section 6 we
prove
an auxiliary result which is then used in section 7 to extend
Theorem B to
the groups $R_{F'/F}{\roman{GL}}_n$ and $R_{F'/F}
{\roman{GSp}}_{2n}$.
\par
Our motivation for the results proved in this paper comes from the
fact
that they make it possible to reformulate in many cases the
conjecture in
[LR] on the reduction of Shimura varieties. Whereas in loc.cit.\ the
concept of {\it admissible morphisms of Galois gerbs} was defined
using
the Bruhat-Tits building, it is possible to replace that condition by
imposing on the corresponding element $b\in B(G)$ that it lie in the
subset $B(G, \mu)$. Here $B(G, \mu)$ is the finite subset of
$B(G)$
defined by the group-theoretic version of Mazur's theorem [K II],
[RR].
The possibility of such a reformulation is implicitly behind the
considerations in section 6 of [K II].
\par
When we presented these results at the Raynaud conference in
Paris,
Fontaine pointed out to us that Theorem A was known to him
earlier in a
different guise (in the case of ${\roman{GL}}_n$). Namely, he had
established the existence of a weakly admissible filtration of type
$\mu$
on the isocrystal $N$, provided that $\mu\geq \boldsymbol\nu(N,F)$. From this
the
existence of the lattice $M$ follows by appealing to the theorem of
Laffaille.
\par
M.~R. wishes to thank the department of mathematics
of the University of Chicago for its hospitality during his visit in the
fall of 2000, when the results of this paper were obtained. He also
thanks
the department of mathematics of the University of Minnesota for
the
possibility of presenting these results.

\medskip\noindent
{\bf Notation.}
\par\noindent
$F$ a finite extension of ${\bold Q}_p$
\par\noindent
$L$ the completion of the maximal unramified extension of $F$
\par\noindent
${O_F}$ resp.\ $O_L$ the rings of integers of $F$ resp.\ $L$
\par\noindent
$\pi\in O_F$ a uniformizer
\par\noindent
${\bold F}$ the residue field of $O_L$
\par\noindent
$\sigma$ the relative Frobenius automorphism of $L/F$.
\medskip\noindent
We follow the tradition of denoting a $\sigma$-linear automorphism
of an
$L$-vector space by $F$ (from ``Frobenius''); there should be no
danger of
confusing this with the notation for the ground field $F$.
\bigskip
\noindent {\bf 1. The result for ${\roman{GL}}_n$}
\medskip\noindent
Let $(N,F)$ be an $F$-isocrystal, i.e.\ a finite-dimensional $L$-
vector
space with a $\sigma$-linear bijective endomorphism. Let $n$
denote the
dimension of $N$. To the $F$-isocrystal $(N,F)$ is associated its
slope
vector $\boldsymbol\nu=\boldsymbol\nu(F)=(\nu_1,\ldots,\nu_n)\in
({\bold Q}^n)_+$. Here $({\bold Q}^n)_+=\{ (\nu_1,\ldots, \nu_n)\in
{\bold
Q}^n;\ \nu_1\geq \nu_2\geq\ldots\geq \nu_n\}$.
\par
Fix an integer $r$ with $0\leq r\leq n$. We call the $F$-isocrystal
$(N,F)$ {\it minuscule of weight} $r$ if the slope vector $\boldsymbol\nu
= \boldsymbol\nu(F)$ of $(N,F)$ satisfies the following
condition $$0\leq \nu_n\leq\ldots\leq\nu_1\leq 1\ \ ,\ \
\sum_{i=1}^n\nu_i=r\ \ .\leqno(1.1)$$ An equivalent condition is the
following. Let $\omega_r$ be the vector $(1,\ldots,1,0,\ldots,0)$,
where 1
is repeated $r$ times and 0 is repeated $n-r$ times. On $({\bold
Q}^n)_+$
we have the usual {\it dominance order,} for which $\nu\leq\mu$ if
and
only if $$\leqalignno{ \nu_1&\leq\mu_1&(1.2)\cr \nu_1+\nu_2 &
\leq\mu_1+\mu_2\cr & \vdots\cr \nu_1+\ldots+ \nu_{n-1}&\leq
\mu_1+\ldots+\mu_{n-1}\cr
\nu_1+\ldots+\nu_n&=\mu_1+\ldots+\mu_n \ \ .\cr
} $$ Then it is easy to see that the condition (1.1) is equivalent to
the
condition $$\boldsymbol\nu(F)\leq \omega_r\ \ .\leqno(1.3)$$ Let
$\bar
I\subset {\bold Z}/n{\bold Z}$ be a non-empty subset and let
$I\subset{\bold Z}$ the inverse image of $\bar I$ under the canonical
surjection ${\bold Z}\to {\bold Z}/n{\bold Z}$. A {\it periodic lattice
chain of type} $\bar I$ in the $L$-vector space $N$ is a set ${\bold
M}$
of $O_L$-lattices $M_i$ $(i\in I)$ for which $$\leqalignno{ & \hbox{if
$i<j$ in $I$, then $M_i\subset M_j$ with length $(M_j/M_i)=j-i$} &
(1.4.1)\cr & M_{i+n}=\pi^{-1}M_i\ \ . & (1.4.2)\cr}$$ In case $\bar
I={\bold Z}/n{\bold Z}$ we also speak of a {\it full periodic lattice
chain.} If $\bar I$ consists of a single element, then a periodic
lattice
chain of type $\bar I$ is simply given by a lattice (namely $M_i$ for
the
unique $i\in I$ with $0\leq i<n$). We denote by $X_{\bar I}$ the set
of
periodic lattice chains of type $\bar I$.
\par
We now fix an isocrystal $(N,F)$ of dimension $n$. We denote by
$X(\omega_r, F)_{\bar I}$ the set of periodic lattice chains of type
$\bar
I$ in $N$ which satisfy the following condition, $$\hbox{for all $i\in I$
we have $\pi M_i\subset FM_i\subset M_i$ and
${\roman{dim}}_{\bold F}\,
M_i/FM_i=r$.}\leqno(1.5)$$ An equivalent condition is the following.
By
the elementary divisor theorem we can associate to any pair of
$O_L$-lattices $M, M'$ in $N$ their relative position
${\roman{inv}}(M,M')=\boldsymbol\mu=(\mu_1,\ldots,\mu_n)\in
({\bold
Z}^n)_+$. Here $({\bold Z}^n)_+={\bold Z}^n\cap ({\bold Q})_+$.
Then the
condition (1.5) is equivalent to $$\hbox{for all $i\in I$ we have
${\roman{inv}}(M_i, FM_i)=\omega_r$}.\leqno(1.6)$$ Also it is clear
that
it suffices to check the conditions (1.5) and (1.6) on a set of
representatives of $I$ mod $n$.
\par
For a non-empty subset $\bar{J}$ of $\bar I$ there is an obvious
forgetful
map $$X(\omega_r, F)_{\bar I}\longrightarrow X(\omega_r,
F)_{\bar{J}}\ \
.\leqno(1.7)$$ We may now formulate the main result of this
section.
\medskip\noindent
{\bf Proposition 1.1.} {\it (i) For any non-empty $\bar I\subset {\bold
Z}/n{\bold Z}$ we have $$X(\omega_r,F)_{\bar I}\neq \emptyset\
\hbox{if
and only if $F$ is minuscule of weight $r$.}$$ (ii) For any non-empty
subsets $\bar I$ and $\bar {J}$ of ${\bold Z}/n{\bold Z}$ with
$\bar{J}\subset \bar I$, the natural map \rom{(1.7)} is surjective. }
\medskip\noindent
To prove this proposition we make the following preliminary
remarks.
\medskip\noindent a) Let $\bar I$ consist of a single element. Then
the statement ``$X(\omega_r, F)_I\neq \emptyset \Longrightarrow
F$ is
minuscule of weight $r$'' is exactly the content of Mazur's
theorem that
the Hodge polygon of an $F$-crystal lies below the Newton polygon
of its
associated $F$-isocrystal and has the same endpoint (use the
reformulations (1.3) resp.\ (1.6) of the relevant conditions).
\medskip\noindent b) If $X(\omega_r, F)_{\bar I}\neq\emptyset$ and
$\bar{J}$ is a non-empty subset of $\bar I$, then obviously
$X(\omega_r,F)_{\bar{J}}\neq\emptyset$.
\medskip
Taking into account a) and b) we see that (i) in Proposition 1.1
follows
from (ii) and the following lemma.
\medskip\noindent
{\bf Lemma 1.2.} {\it Let $F$ be minuscule of weight $r$. Then there
exists an $O_L$-lattice $M$ in $N$ with $${\roman{inv}}
(M,FM)=\omega_r\ \
.$$ }
\medskip\noindent
{\bf Proof.} Let us first assume that $F$ is isoclinic, i.e.\ $\boldsymbol\nu(F)=
(\nu,\ldots,\nu)$ with $0\leq \nu\leq 1$ and $n\nu=r$. In this case
the
$F$-isocrystal is uniquely determined up to isomorphism and there
exists a
basis $e_1,\ldots,e_n$ of $N$ such that $$Fe_1=e_2,\
Fe_2=e_3,\ldots,
Fe_{n-1}=e_n,\ Fe_n=\pi^re_1\ \ . \leqno(1.8)$$ Then the following
lattice
is as required, $$M=O_L\cdot\pi^{r-1}\cdot e_1\oplus O_L\pi^{r-
2}\cdot
e_2\oplus\ldots\oplus O_L\pi\cdot e_{r-1}\oplus O_L\cdot
e_r\oplus\ldots\oplus O_L\cdot e_{n-1}\oplus O_L\cdot e_n\ .$$
The general
case follows since the isocrystal $(N,F)$ is the direct sum of
isoclinic
isocrystals $(N_i, F_i)$ $(i=1,\ldots, s)$ which are minuscule of
weight
$r_i$, with $\sum_{i=1}^sr_i=r$.\endbeweis
\medskip
To prove (ii) of Proposition 1.1, we may assume that $\bar I={\bold
Z}/n{\bold Z}$. Hence starting from $\bar J$ we may enlarge $\bar
J$ by
one element at a time. We are then reduced to proving the following
statement.
\medskip\noindent
{\bf Lemma 1.3.} {\it Consider $O_L$-lattices $M,M'$ such that
$$M\mathop{\supset}\limits_{\neq} M'\supset \pi M\ \ ,$$ with
${\roman{inv}}(M,FM)=\omega_r$, ${\roman{inv}}(M',
FM')=\omega_r$. Then
there exists an $O_L$-lattice $\tilde M$ such that $$M\supset \tilde
M\supset M'\ \ ,$$ with ${\roman{dim}}_{\bold F}\tilde M/M'=1$ and
${\roman{inv}}(\tilde M, F\tilde M)=\omega_r$.}
\medskip\noindent
{\bf Proof.} We introduce the $\sigma^{-1}$-linear operator $V$
defined by
the identity $$VF=FV=\pi\ \ .\leqno(1.9)$$ Then, since $F$ is
minuscule of
weight $r$, the condition ${\roman{inv}}(M,FM)=\omega_r$ on a
lattice $M$
is equivalent to the condition $$FM\subset M\ \ \hbox{and}\ \
VM\subset M\
\ .\leqno(1.10)$$ Consider the ${\bold F}$-vector space $W=M/M'$
with the
induced $\sigma^{\pm 1}$-linear operators $\overline F, \overline V$
which
satisfy $\overline F\, \overline V=\overline V\, \overline F=0$. By the
previous remarks it suffices to find a line $\ell$ in $W$ which is
stable
under $\overline F$ and $\overline V$. We distinguish cases.
\medskip\noindent
{\bf Case 1.} $\overline F$ is bijective.
\par\noindent
In this case there exists an ${\bold F}$-basis of $W$ consisting of
$\overline F$-invariant vectors (Dieudonn\'e). Let $\ell$ be the line
generated by one of these basis vectors. Since $\overline V=0$ in
this
case, this line is stable under $\overline F$ and $\overline V$.
\medskip\noindent
{\bf Case 2.} Ker $\overline F\neq(0)$.
\par\noindent
The map $\overline V$ induces a map from  Ker $\overline F$ to
itself. If this induced map fails
to be bijective, we take $\ell$ to be any line in its kernel. If the
induced map is
bijective, so that there exists a basis of  Ker $\overline F$
consisting of $\overline V$-invariant
vectors, then we take
$\ell$ to be the line generated by one of the basis vectors.
\endbeweis

\bigskip\noindent \noindent {\bf 2. The result for
${\roman{GSp}}_{2n}$}
\medskip\noindent
Let $(N, \langle\ ,\ \rangle)$ be a symplectic vector space of
dimension
$2n$ over $L$. Let $F$ be a $\sigma$-linear bijective
endomorphism of $N$
satisfying $$\langle Fx,Fy\rangle =c\cdot \langle
x,y\rangle^{\sigma}\ \
,\ \ x,y\in N\leqno(2.1)$$ for some fixed $c\in L^{\times}$. We call
$(N,
\langle\ ,\ \rangle, F)$ a {\it symplectic} $F${\it -isocrystal.} The
slope vector $\boldsymbol\nu( F)$ of the isocrystal $(N,F)$ then
satisfies $$\nu_1+\nu_{2n}=\nu_2+\nu_{2n-1}=\ldots
=\nu_n+\nu_{n+1}=d\ \
,\leqno(2.2)$$ where $d={\roman{val}}(c)$ is the $\pi$-adic valuation
of
$c$. We call the symplectic $F$-isocrystal {\it minuscule of
weight} $r$
for some $r$ with $0\leq r\leq 2n$ if the underlying $F$-isocrystal is
minuscule of weight $r$ in the sense of (1.1). Note that only $r=0$,
$r=n$
and $r=2n$ are possible and that then $d=0,1$ or 2 respectively.
\par
Let $\bar I$ be a non-empty {\it symmetric} subset of ${\bold
Z}/2n{\bold
Z}$, i.e., invariant under multiplication by $-1$. Let $I$ be the
inverse
image of $\bar I$ under the surjection ${\bold Z}\to {\bold Z}/2n{\bold
Z}$. A periodic lattice chain $M_i$ $(i\in I)$ of type $\bar I$ is called
{\it selfdual} if  there exists $d\in{\bold Z}$ such that
$$M_i^{\perp}=M_{-i+d\cdot 2n}\ \ ,\ \ i\in I\ \ .\leqno(2.3)$$ Here for
any $O_L$-lattice $M$ in $N$ we put $$M^{\perp}=\{ x\in N;\
\langle x,
M\rangle\subset O_L\}\ \ .\leqno(2.4)$$ We denote by $X^G_{\bar
I}$ the
set of selfdual periodic lattice chains of type $\bar I$ in $N$. Let
now
$(N, \langle\ ,\ \rangle, F)$ be a symplectic $F$-isocrystal. For a
non-empty symmetric subset $\bar I$ in ${\bold Z}/2n{\bold Z}$, let
$X^G(\omega_r, F)_{\bar I}$ denote the set of periodic selfdual
lattice
chains of type $\bar I$ in $N$ which lie in the set $X(\omega_r,
F)_{\bar
I}$ in the sense of (1.5) for ${\roman{GL}}_{2n}$.
\medskip
The following result is the analogue of Proposition 1.1 in the
present context.
\medskip\noindent
{\bf Proposition 2.1.} {\it (i) For any non-empty symmetric subset
$\bar
I$ of ${\bold Z}/2n{\bold Z}$ we have $$\eqalign{ X^G(\omega_r,
F)_{\bar
I}\neq \emptyset\ &\hbox{if and only if the symplectic $F$-isocrystal
$(N,\langle\ ,\ \rangle, F)$ is minuscule} \cr &\hbox{of weight
$r$.}\cr}$$ (ii) For any non-empty symmetric subsets $\bar I$ and
$\bar J$
of ${\bold Z}/2n{\bold Z}$ with $\bar J\subset\bar I$, the natural map
$$X^G(\omega_r, F)_{\bar I}\longrightarrow X^G(\omega_r, F)_{\bar
J}$$ is
surjective. }
\medskip\noindent
Again, by Mazur's theorem, we infer that if $X^G(\omega_r, F)_{\bar
I}\neq\emptyset$, then $F$ is minuscule of weight $r$ (in particular
$r=0$, or $n$, or $2n$). Conversely, assume that $F$ is minuscule
of
weight $r$. If $r=0$, then $N$ admits a symplectic basis of $F$-
invariant
vectors (Dieudonn\'e), hence defines an $F$-form $(N_0, \langle\ ,\
\rangle_0)$ of $(N,\langle\ ,\ \rangle)$. Any self-dual $O_F$-lattice
in
$N_0$ defines an element of $X^G(\omega_0, F)_{\bar I}$, where
$\bar I=\{
0\}$. Furthermore, the assertion (ii) of Proposition 2.1 just amounts
to
the fact that any selfdual periodic lattice chain may be completed
to a
full selfdual periodic lattice chain. This is well-known, comp.\ [KR],
section 10. The case $r=2n$ reduces to the previous one by
replacing $F$
by $\pi^{-1}\cdot F$. Hence from now on we may assume that $F$
is
minuscule of weight $n$.
\medskip\noindent
{\bf Lemma 2.2.} {\it Let $(N,\langle\ ,\ \rangle, F)$ be a symplectic
$F$-isocrystal which is minuscule of weight $n$. Then there exists
a
selfdual $O_F$-lattice $M$ such that $$M\supset FM\supset\pi M\ \
\hbox{and}\ \ (FM)^{\perp}=\pi^{-1}FM\ \ .$$ In other words $M\in
X^G(\omega_n, F)_{\bar I}$ with $\bar I=\{ 0\}$. }
\medskip\noindent
{\bf Proof.} By hypothesis $0\leq \nu_{2n}\leq \nu_{2n-
1}\leq\ldots\leq
\nu_1\leq 1$ and
$$\nu_1+\nu_{2n}=\nu_2+\nu_{2n-1}=\ldots=\nu_n+\nu_{n+1}=1\ \
.$$ From the
slope decomposition of $N$ we deduce a direct sum
decomposition $$N=
N'\oplus \tilde N\oplus N''\ \ ,$$ where $N'$ resp.\ $N''$ includes all
slope components of slope $<1/2$ resp.\ $>1/2$ and where $\tilde
N$ is the
sum of all slope components of slope $1/2$. Then $N'$ and $N''$
are
totally isotropic subspaces which are in duality by $\langle\ ,\
\rangle$
and $\tilde N$ is orthogonal to $(N'\oplus N'')$. A selfdual lattice in
$N'\oplus N''$ may be obtained by taking any $O_F$-lattice $M'$ in
$N'$
and then forming $M'\oplus M''$ where $$M''= M^{'\perp}=\{ x\in N'';\
\langle x,M'\rangle \subset O_L\}\ \ .$$ Using the result of section 1
for
${\roman{GL}}_{n'}$, where $n'={\roman{dim}}\ N'$, we are reduced
to
considering $\tilde N$, i.e., we may assume from the start that all
slopes
of $N$ are equal to $1/2$. In this case the symplectic $F$-
isocrystal is
uniquely determined up to isomorphism and there exists a basis of
$N$ such
that $$\leqalignno{ & Fe_i=-e_{2n-i+1},\ Fe_{2n-i+1}=\pi e_i\ \
i=1,\ldots,n\ \hbox{and} \cr & \langle e_i, e_j\rangle=0, \ \langle
e_{2n-i+1}, e_{2n-j+1}\rangle =0,\ \langle e_i, e_{2n-j+1}\rangle
=\delta_{ij},\ \ i,j= 1,\ldots,n\ \ . \cr}$$ Then the $O_L$-lattice $M$
generated by $e_1,\ldots, e_{2n}$ satisfies the required conditions.
\endbeweis
\medskip
To complete the proof of Proposition 2.1, it suffices now to prove
assertion (ii) in the case where $\bar I= {\bold Z}/2n{\bold Z}$.
Enlarging $\bar J$ one step at a time we then reduce to the case in
which
$\bar J\subset \bar I$ is as in [KR], 10.2. In other words, we fix
$k\in
J$ such that $k+1\not\in J$ and obtain $\bar I$ by adding to $\bar
J$ one
or two elements, namely the class(es) of $k+1$ and $-(k+1)$
modulo
$2n{\bold Z}$.
\medskip\noindent
Let $\ell$ be the smallest integer in $J$ such that $\ell >k$; thus
$k<\ell\leq k+2n$. Since $\bar J\subset\bar I$, there is a natural
map
$$f:X^G_{\bar I}\longrightarrow X^G_{\bar J}\ \ .\leqno(2.5)$$ We are
interested in the fiber $f^{-1}({\bold M})$ over an element ${\bold
M}=(M_i)_{i\in J}$ of $X^G_{\bar J}$. We associate to an element
${\tilde{\bold M}}=(M_i)_{i\in I}$ of $f^{-1}({\bold M})$ the lattice $M:=
M_{k+1}$. Clearly this lattice satisfies $$\leqalignno{ & M_k\subset
M\subset M_{\ell} & (2.6.1) \cr & {\roman{dim}}_{\bold F} M/M_k=1\
\ . &
(2.6.2) \cr}$$
\medskip\noindent
{\bf Lemma 2.3.} {\it The map $\tilde{\bold M}\mapsto M$ is a
bijection
from $f^{-1}( {\bold M})$ to the set of lattices $M$ in $N$ satisfying
\rom{(2.6.1)} and \rom{(2.6.2)}. }
\medskip\noindent
{\bf Proof.} One way to prove the lemma would be to appeal to general
results of
Bruhat-Tits. In the special case at hand it is also easy to give an
elementary proof, as we
now do.

The map is obviously injective since $\tilde{\bold M}$
contains with $M$ also $M^{\perp}$ and all multiples of these two
lattices. To prove surjectivity we start with $M$ satisfying (2.6.1)
and
(2.6.2) and have to construct $\tilde{\bold M}\in f^{-1}({\bold M})$
which
gives $M$. We imitate the proof of [KR], Lemma 10.3.
\par
Let $\overline P=\bar J\cup \{ \overline k+\overline 1\}$ (for $m\in
{\bold Z}$, we write $\overline m$ for its class modulo $2n$). Let
$\overline Q =-\overline P$. Then $\bar I=\overline P\cup\overline
Q$; and
for the inverse images $P$ and $Q$ of $\overline P$ and $\overline
Q$ in
${\bold Z}$, we have $P=-Q$.
\par
There is a unique periodic lattice chain ${\bold X}$ of type $\overline
P$
such that $X_{k+1}=M$ and $X_j=M_j$ for $j\in J$. Let $d\in{\bold
Z}$ be
the unique integer such that $M_j^{\perp}=M_{-j+d\cdot 2n}$ for
$j\in J$.
There is a unique periodic lattice chain ${\bold Y}$ of type
$\overline Q$
such that $Y_{-(k+1)}=\pi^dM^{\perp}$ and $Y_j=M_j$ for $j\in J$.
\par
We claim that $$p\in P\ ,\  q\in Q\  ,\ p\leq q\Longrightarrow
X_p\subset
Y_q\ \ .\leqno(2.7)$$ This is obvious if there exists $j\in J$ such that
$p\leq j\leq q$, so we now assume the contrary. It is harmless to
suppose
that $p=k+1$. Then necessarily $q=\ell -1$ and $\overline{\ell}=-
\overline
k$. Consider the ${\bold F}$-vector space $V=M_{\ell}/M_k$. Then
the
lattices $X_p$ resp.\ $Y_q$ correspond to subspaces $U_1$ resp.\
$U_2$ of
$V$, where $U_1$ is of dimension one and $U_2$ is of
codimension one. We
have to show that $U_1\subset U_2$. But on $V$ we have the
symplectic form
defined by the fact that $M_{\ell}=\pi^r\cdot M_k^{\perp}$, where $r$
is
defined by $\ell =-k+r\cdot 2n$. Furthermore, we have $Y_{\ell -1}=
\pi^r\cdot X^{\perp}_{k+1}$. Equivalently, we have
$U_2=U_1^{\perp}$ for
the symplectic form on $V$. The claim now follows from the fact
that any
line in a symplectic vector space is isotropic.
\par
We also claim that $$p\in P\ ,\ q\in Q\ \ q<p\Longrightarrow
Y_q\subset
X_p\ \ .\leqno(2.8)$$ This is clear since there always exists $j\in J$
such that $q\leq j\leq p$.
\par
Now suppose $p\in P$, $q\in Q$ and $p=q$. Then from (2.7) we
have
$X_p\subset Y_q$. But both are lattices which contain $M_{k-
r\cdot 2n}$
for sufficiently large $r$ and with the same index, hence
$X_p=X_q$. Thus,
without ambiguity, we may define the periodic lattice chain
$\tilde{\bold
M}=(M_i)_{i\in I}$ of type $\bar I$ by putting $M_i=X_i$ if $i\in P$
and
$M_i=Y_i$ if $i\in Q$. It is obvious that this is indeed a selfdual
lattice chain contained in $f^{-1}({\bold M})$ and that $\tilde{\bold
M}\mapsto M$.\endbeweis
\medskip\noindent
Using this lemma, the surjectivity assertion (ii) in Proposition 2.1 is
reduced to the corresponding statement for ${\roman{GL}}_{2n}$, which is
Lemma 1.3.\endbeweis
\bigskip \noindent
{\bf 3. The general problem}
\medskip\noindent
Let $G$ be a connected reductive group over $F$. For simplicity we assume
that $G$ splits over $L$. (The problem addressed in this section can be
formulated without this hypothesis, but then becomes more technical and
even more speculative). Let $\tilde T$ be a maximal split torus over $L$.
Let ${\Cal B}= {\Cal B}(G_{ad}, L)$ be the Bruhat-Tits building of the
adjoint group over $L$. To $\tilde T$ corresponds an apartment in ${\Cal
B}$. Let $\tilde K_0$ be an Iwahori subgroup of $G(L)$ corresponding to an
alcove in the apartment for $\tilde T$. Let $\tilde W$ be the {\it Iwahori
Weyl group of} $\tilde T$, $$\tilde W=\tilde N(L)/\tilde T(L)_1\ \
.\leqno(3.1)$$ Here $\tilde N$ denotes the normalizer of $\tilde T$ and
$\tilde T(L)_1$ the maximal bounded subgroup of $\tilde T(L)$. Then
$\tilde T(L)_1=\tilde T(L)\cap \tilde K_0$. Let $\tilde K$ be the
parahoric subgroup of $G(L)$ corresponding to a facet of the base alcove.
 Let $$\tilde W^{\tilde K}=\tilde N(L)\cap \tilde K/\tilde T(L)\cap \tilde K_0\
\ .\leqno(3.2)$$ Then there is a canonical bijection $$\tilde K\setminus
G(L)/\tilde K= \tilde W^{\tilde K}\setminus \tilde W/\tilde W^{\tilde K}\
\ .\leqno(3.3)$$ We therefore obtain a succession of maps whose
composition will be denoted by inv, \par $$
%\gather
{\roman{inv}}: G(L)/\tilde K\times G(L)/\tilde K{\to} G(L)\setminus
(G(L)/\tilde K\times G(L)/\tilde K) {=}\tilde K\setminus G(L)/\tilde K{=}
\tilde W^{\tilde K}\setminus \tilde W/\tilde W^{\tilde K}.\tag3.4
%\endgather
$$ We now fix a conjugacy class of minuscule one-parameter subgroups $\mu$
of $G$ defined over $L$. We may assume that $\mu$ factors through $\tilde
T$ and determines an orbit in $X_*(\tilde T)$ under the conjugation action
of the finite Weyl group $W=\tilde N(L)/\tilde T(L)$. Let
${\roman{Adm}}(\mu)\subset \tilde W$ be the {\it admissible subset
corresponding to} $\mu$ ([KR], Introduction), $${\roman{Adm}}(\mu)=\{ w\in
\tilde W;\ w\leq t_{\mu'},\ \hbox{some}\ \mu'\}\ \ .\leqno(3.5)$$ Here
$\mu'$ denotes an element of the $W$-orbit in $X_*(\tilde T)$ defined by
$\mu$, and $t_{\mu'}$ the corresponding element of $\tilde W$. In (3.5)
appears the Bruhat order on $\tilde W$ determined by the base alcove. We
denote by $${\roman{Adm}}_{\tilde K}(\mu)\subset \tilde W^{\tilde
K}\setminus \tilde W/\tilde W^{\tilde K}\leqno(3.6)$$ the image of
${\roman{Adm}}(\mu)$ under the natural projection. It is independent of
the choice of $\tilde K_0$ contained in $\tilde K$. We will assume that
$\tilde K$ is $\sigma$-invariant, or equivalently that the corresponding
facet in the building is $\sigma$-invariant. Then $K=\tilde
K^{\langle\sigma\rangle}$ is a parahoric subgroup of $G(F)$. We note that,
conversely, $K$ determines $\tilde K$ and the corresponding
$\sigma$-invariant facet in ${\Cal B}$.
\medskip
Our final choice is an element $b\in G(L)$. We then define $$X(\mu,b)_K=\{
g\in G(L)/\tilde K;\ {\roman{inv}}(g,b\sigma(g))\in{\roman{Adm}}_{\tilde
K}(\mu)\}\ .\leqno(3.7)$$ Let $\tilde K'$ be a $\sigma$-invariant
parahoric subgroup of $G(L)$ containing $\tilde K$. Then $K'=\tilde
K^{'\langle\sigma\rangle}$ is a parahoric subgroup of $G(F)$ and there is
a canonical projection map $$X(\mu, b)_K\longrightarrow X(\mu, b)_{K'}\ \
.\leqno(3.8)$$ Let $B(G)$ be the set of $\sigma$-conjugacy classes in
$G(L)$ and let $[b]\in B(G)$ be the $\sigma$-conjugacy class of $b$. We
denote by $B(G,\mu)$ the finite subset of $B(G)$ defined by the group
theoretic version of Mazur's theorem ([K II], \S 6).
\medskip\noindent
{\bf Conjecture 3.1.} {\it (i) For any parahoric subgroup $K$ of $G(F)$ we
have $$X(\mu, b)_K\neq\emptyset\Longleftrightarrow [b]\in B(G,\mu)\ \ .$$
(ii) For any pair of parahoric subgroups $K\subset K'$ of $G(F)$, the map
(3.8) is surjective. }
\medskip\noindent
It is not clear whether the hypothesis that $\mu$ is minuscule is indeed
necessary for the statements in this Conjecture. \par Let
$G={\roman{GL}}_n$. A conjugacy class of minuscule one-parameter subgroups
of $G$ is of the form $\mu=\omega_r+k\cdot\omega_n= \omega_r+k\cdot {\bold
1}$ for a unique $r$ with $0\leq r<n$ and some $k\in {\bold Z}$. Here
${\bold 1}=(1,1,\ldots, 1)$. The validity of Conjecture 3.1 is unchanged
if $\mu$ is replaced by $\omega_r$, so we assume this now.
\par
The conjugacy classes of parahoric subgroups correspond in a one-to-one
way to the set of non-empty subsets $\bar I\subset {\bold Z}/n{\bold Z}$
and the corresponding coset space may be identified with the space
$X_{\bar I}$ of periodic lattice chains of type $\bar I$. Let ${\bold
M}=(M_i)_{i\in I}$ and ${\bold M}'=(M'_i)_{i\in I}$ be two elements of
$X_{\bar I}$. Then $${\roman{inv}}({\bold M}, {\bold M}')\in
{\roman{Adm}}_{\tilde K_I}(\mu)\Leftrightarrow M_i\supset M'_i\supset \pi
M_i\ \hbox{and}\ {\roman{dim}}_{\bold F}M_i/M'_i=r,\ \forall i\in I\
.\leqno(3.9)$$ Indeed, this follows from the identification of
${\roman{Adm}}_{\tilde K_I}(\mu)$ with the $\mu$-permissible set inside
$\tilde W^{K_I}\setminus \tilde W/\tilde W^{K_I}$, [KR], [HN]. In fact,
for any dominant coweight $\mu$ we have (comp.\ [HN], 9.7)
$${\roman{inv}}({\bold M}, {\bold M}')\in {\roman{Perm}}_{\tilde
K_I}(\mu)\Longleftrightarrow {\roman{inv}}(M_i, M'_i)\leq \mu\ \ ,\ \
\forall i\in I\ \ .\leqno(3.10)$$ If $\mu$ is minuscule, the inequality on
the right hand side is necessarily an equality which yields the condition
appearing in (3.9). These remarks imply that the results of section 1
prove Conjecture 3.1 in the case of ${\roman{GL}}_n$.
\par
Similarly, the results of section 2 prove Conjecture 3.1 in the case of
$G={\roman{GSp}}_{2n}$. In fact, in this case the $\mu$-admissible set is
the intersection of the $\mu$-permissible set for ${\roman{GL}}_{2n}$ with
the extended affine Weyl group of ${\roman{GSp}}_{2n}$, cf.\ [KR], see
also [HN], Prop.\ 9.7.

\bigskip\noindent
{\bf 4. A converse to Mazur's inequality}
\medskip\noindent
In this section we let $G={\roman{GL}}_n$ or
$G={\roman{GSp}}_{2n}$. Our
aim is to prove a converse to Mazur's theorem, strengthening for
these
groups Prop.\ 4.2.\ of [R]. Much of our argument remains valid
for an
arbitrary split group with simply connected derived group.
\par
We start with a lemma which is the group-theoretic interpretation of
the
first half of the proof of Lemma 1.2. Let $A$ be a maximal split
torus in
$G$. We denote by $\pi_1(G)$ the algebraic fundamental group of
$G$. Since
$G_{\roman{der}}$ is simply connected, $\pi_1(G)$ is the factor
group of
$X_*(A)$ by the lattice generated by the coroots, and is a free
abelian
group. We denote by $$\kappa_G:G(L)\longrightarrow
\pi_1(G)\leqno(4.1)$$
the homomorphism introduced in [K II]. We denote by $\tilde K= G(
O_L)$
the special maximal bounded subgroup determined by a Chevalley
form of $G$
adapted to $A$.
\medskip\noindent
{\bf Proposition 4.1.} {\it Let $g\in G(L)$ and let $b\in G(L)$ be a
basic
element. Then the $\sigma$-conjugacy class of $b$ meets $\tilde
Kg\tilde
K$ if and only if $\kappa_G(g)=\kappa_G(b)$.
\par\noindent
  }
\medskip\noindent
{\bf Proof.} One direction is trivial, since $\kappa_G(\tilde K)=\{ 0\}$
and since $\sigma$-conjugate elements have identical images
under
$\kappa_G$. For the converse direction we may use the Cartan
decomposition
of $G(L)$ to assume that $g\in A(L)$ and even $g=a\in A(F)$. Let
$w\in W$
be an elliptic element, i.e.\ $X_*(A)_{\bold
R}^w=X_*(Z_G^{\circ})_{\bold
R}$. Here $Z_G^{\circ}$ denotes the connected center of $G$.
Equivalently,
any $w$-invariant element of $A(F)$ has finite order modulo the
center.
Let $\dot w\in N_G(A)(F)\cap G(O_F)$ be a representative of $w$ in
$G(O_F)$.
We claim that $a\dot w$ is a basic element in $G(L)$.
Once this is
established, we conclude from $\kappa_G(a\dot
w)=\kappa_G(a)=\kappa_G(b)$
that $a\dot w$ and $b$ are $\sigma$-conjugate ([K I], 5.6), which
finishes
the proof since $a\dot w\in\tilde K a\tilde K$.
\par
To see that $a\dot w$ is basic it suffices to show that its norm
under a
sufficiently large finite extension $F'$ of $F$ contained in $L$ is
central ([K I], 4.3.). Since $a\dot w\in G(F)$, we have to see that a
sufficiently high power of $a\dot w$ is central. But $$(a\dot
w)^r=a\cdot
w(a)\cdot\ldots\cdot w^{r-1}(a)\cdot \dot w^r\ \ .\leqno(4.2)$$ If $r$ is
divisible by the order of $w$ in $W$ we have that
$aw(a)\cdot\ldots\cdot
w^{r-1}(a)$ is $w$-invariant and hence is of finite order modulo the
center. The same applies to $\dot w^r$ and hence our claim is
proved.\endbeweis

\medskip\noindent In the sequel we fix a Borel subgroup $B=A U$.
We denote
by $X_*(A)_{\roman{dom}}$ resp.\ $X_*(A)_{{\bold Q},
{\roman{dom}}}$ the
set of dominant elements in $X_*(A)$ resp.\ $X_*(A)\otimes{\bold
Q}$.
Recall ([K II], 4.2) that to $b\in G(L)$ is associated its Newton point
$\overline{\nu} (b)\in X_*(A)_{{\bold Q},{\roman{dom}}}$. We will
denote
by $\leq$ the usual partial order on $X_*(A)_{{\bold
Q},{\roman{dom}}}$,
i.e.\ $\nu\leq \nu'$ iff $\nu'-\nu$ is a non-negative linear combination
of positive coroots. Note that, since the derived group of $G$ is
simply
connected, the partial order induced on $X_*(A)_{\roman{dom}}$ is
that
denoted in [R] by $\mathop{\leq}\limits^!$, i.e.\
$\nu\mathop{\leq}\limits^!\nu'$ iff $\nu'-\nu$ is a non-negative {\it
integral} linear combination of positive coroots.
\medskip\noindent
{\bf Corollary 4.2.} {\it Let $b\in G(L)$ be basic with associated
Newton
point $\overline{\nu}=\overline{\nu}(b)\in X_*(A)_{{\bold
Q},{\roman{dom}}}$. Let $\mu\in X_*(A)_{\roman{dom}}$ with
$\overline{\nu}\leq \mu$. Then there exists $h\in G(L)/\tilde K$ with
$${\roman{inv}}(h, b\sigma(h))=\mu\ \ .$$ }
\par\noindent {\bf Proof.} Let $g=\pi^{\mu}\in A(F)$. Then
$\kappa_G(g)=\kappa_G(b)$ and applying the previous proposition
we find
$h\in G(L)$ with $h^{-1}b\sigma (h)\in \tilde K\pi^{\mu}\tilde K$, as
desired.\endbeweis
\bigskip\noindent {\bf Remark 4.3.} In the case of
${\roman{GL}}_n$ the previous construction can be made totally
explicit.
In this case $\pi_1(G)={\bold Z}$ and any basic $b\in G(L)$ with
$\kappa_G(b)=r\in{\bold Z}$ is $\sigma$-conjugate to the element
$F$
described by (1.8). Let $\mu\in ({\bold Z}^n)_+$ with
$\sum_{i=1}^n\mu_i=r$. Then the lattice $M$ spanned by the
vectors
$$\pi^{\sum^n_{i=2}\mu_i}e_1,\pi^{\sum^n_{i=3}\mu_i}e_2,\ldots,
\pi^{\mu_n}e_{n-1},e_n\leqno(4.3)$$ satisfies
${\roman{inv}}(M,FM)=\mu$.

\medskip
Let now $P=M N$ be a parabolic subgroup containing $B$, where
$M$ is the
unique Levi subgroup containing $A$. We sometimes consider
$M$ as a factor
group of $P$. For $\mu\in X_*(A)$ we denote by $M(\mu)$ the
image of
$\tilde K\pi^{\mu}\tilde K\cap P(L)$ in $M(L)$.

\medskip\noindent
{\bf Lemma 4.4.} {\it Let $b\in M(L)$ and let $\mu\in X_*(A)$. Then
the
$\sigma$-conjugacy class of $b$ in $G(L)$ meets $\tilde
K\pi^{\mu}\tilde
K$ if and only if the $\sigma$-conjugacy class of $b$ in $M(L)$
meets
$M(\mu)$. }
\medskip\noindent
{\bf Proof.} Assume that the $\sigma$-conjugacy class of $b$
meets $\tilde
K\pi^{\mu}\tilde K$. By the Iwasawa decomposition, there then
exists $p\in
P(L)$ with $pb\sigma(p)^{-1}\in\tilde K\pi^{\mu}\tilde K$. Writing
$p=m.n\in M(L).N(L)$ we conclude that $mb\sigma(m)^{-1}\in
M(\mu)$.
\par
Conversely, assume there exists $m\in M(L)$ with $mb\sigma(m)^{-
1}\in
M(\mu)$. Hence there exists $n\in N(L)$ with $mb\sigma(m)^{-
1}n\in \tilde
K\pi^{\mu}\tilde K$. But by [K II], 3.6, the two elements
$mb\sigma(m)^{-1}n$ and $mb\sigma(m)^{-1}$ are
$\sigma$-conjugate by an
element in $P(L)$. Hence $b$ is $\sigma$-conjugate in $G(L)$ to
an element
in $\tilde K\pi^{\mu}\tilde K$.\endbeweis

\bigskip
Let $\mu\in X_*(A)$ and let $${\Cal P}_{\mu}=\{ \nu\in X_*(A);\
\kappa_G(\nu)=\kappa_G(\mu),\ \nu\in {\roman{Conv}}(W\mu )\}\ \
.\leqno(4.4)$$ Here we have denoted by $\kappa_G:X_*(A)\to
\pi_1(G)$ the
map which sends $\mu$ to $\kappa_G(\pi^{\mu})$. Also
${\roman{Conv}}(W\mu
)$ denotes the convex hull of $W\mu $ in $X_*(A)\otimes {\bold
R}$. Note
that since the derived group of $G$ is simply connected, the first
condition in (4.4) is implied by the second.

\medskip\noindent
{\bf Lemma 4.5.} {\it We have $$\kappa_M(M(\mu))=\kappa_M({\Cal
P}_{\mu})\
\ .$$ }

\medskip\noindent
{\bf Proof.} Let $m\in M(\mu)$ and let us prove that $\kappa_M(m)\in
\kappa_M({\Cal P}_{\mu})$. By the definition of $M(\mu)$ there
exists
$n\in N(L)$ with $mn\in \tilde K\pi^{\mu}\tilde K$. Using the Cartan
decomposition of $M$ we may write $m=k_M\cdot \pi^{\nu}\cdot
k'_M$ with
$k_M, k'_M\in \tilde K_M= M( O_L)$. Then
$\kappa_M(m)=\kappa_M
(\pi^{\nu})$. Now $mn=k_M\pi^{\nu}\cdot n'\cdot k'_M$ with $n'\in
N(L)$.
Hence $\pi^{\nu}n'\in \tilde K\pi^{\mu}\tilde K$. By Satake (comp.\
[R]) this implies $\nu\in {\Cal P}_{\mu}$.
\par
  Conversely, let $\nu\in
{\Cal P}_{\mu}$. Then by [R], Thm.\ 1.1.\ there exists $u\in
U(L)$
such that $\pi^{\nu}u\in \tilde K\pi^{\mu}\tilde K$. Writing $u=u_M.n$
with $u_M\in U(L)\cap M(L)$ and $n\in N(L)$ we have
$\pi^{\nu}\cdot u_M\in
M(\mu)$. But the image of $\nu$ in $\pi_1(M)$ is equal to
$\kappa_M(\pi^{\nu}u_M)$ and hence lies in
$\kappa_M(M(\mu))$.\endbeweis

\medskip\noindent
{\bf Proposition 4.6.} {\it Let $b\in M(L)$ be basic, and let $\mu\in
X_*(A)$. The $\sigma$-conjugacy class of $b$ in $G(L)$ meets
$\tilde
K\pi^{\mu}\tilde K$ if and only if $\kappa_M(b)\in \kappa_M({\Cal
P}_{\mu})$. }
\medskip\noindent
{\bf Proof.} This is a consequence of the results established so far.
Indeed, the $\sigma$-conjugacy class of $b$ in $G(L)$ meets
$\tilde
K\pi^{\mu}\tilde K$ iff the $\sigma$-conjugacy class of $b$ in
$M(L)$
meets $M(\mu)$. Now $M(\mu)$ is a union of $\tilde K_M$-double
cosets.
Applying Proposition 4.1 to each $\tilde K_M$-double coset (with
$M$
instead of $G$), we see that this holds iff $\kappa_M(b)\in
\kappa_M(M(\mu))$. But by the previous lemma we may identify
$\kappa_M(M(\mu))$ and $\kappa_M({\Cal P}_{\mu})$.\endbeweis

\medskip
Recall that $\mu\in X_*(A)$ is called minuscule (in the large sense)
if
$\langle \mu,\alpha\rangle \in\{ 0,\pm 1\}$ for all roots $\alpha$. It is
well-known that $\kappa_G$ induces a bijection (Bourbaki:
Groupes et
Alg\`ebres de Lie, ch. VI, \S 2, ex. 2) $$\{ \mu\in
X_*(A)_{\roman{dom}};\
\mu\ \hbox{minuscule}\} \longrightarrow \pi_1(G)\ \ .\leqno(4.5)$$
\par
Recall our parabolic subgroup $P=M N$. We let $A_M$ be the
maximal split
torus in the center of $M$ and let $X_M=X_*(A_M)\subset X_*(A)$.
Then
$\kappa_M$ induces an injective map
$X_M\hookrightarrow\pi_1(M)$ with
finite cokernel.

\medskip\noindent
{\bf Lemma 4.7.} {\it Let $G={\roman{GL}}_n$.
Let $\mu\in X_*(A)_{\roman{dom}}$ be minuscule and let $x\in
\pi_1(M)$.
The following conditions on $x$ are equivalent:
\par\noindent
\rom{(}i\rom{)} $x\in \kappa_M(W\mu )$
\par\noindent
\rom{(}ii\rom{)} Let $\nu\in X_M\otimes {\bold Q}$ be the unique element
mapping under
$\kappa_M\otimes {\bold Q}$ to $x$. Then $\nu\in
{\roman{Conv}}(W\mu )$.
\par\noindent
If $x$ satisfies these conditions, let $\tilde{\nu}\in X_*(A)$ be the
unique $M$-dominant $M$-minuscule element mapping to $x$, cf.
\rom{(4.5)}. Then
$\tilde\nu\in W\mu $.}
\medskip\noindent
{\bf Proof.} Let $G={\roman{GL}}_n$ and
$M=M_{(m_1,\ldots,m_r)}$. Since
$\mu$ is minuscule, we may write $\mu=k\cdot {\bold
1}+\omega_s$, where
$0\leq s<n$, and $k\in{\bold Z}$ and where we used the notation
${\bold
1}=\omega_n=(1,\ldots, 1)$. Since adding to $\mu$ an element of
$X_G$ does
not affect the assertion of the lemma, we may assume
$\mu=\omega_s$. But
then it is obvious that $$\kappa_M(W\mu )= \{ (s_1,\ldots, s_r);\
0\leq
s_i\leq m_i,\ \Sigma_{i=1}^r s_i=s\}\ \ .\leqno(4.6)$$ Here we have
used
the identification $$\pi_1(M)=
\pi_1({\roman{GL}}_{m_1}\times\ldots\times
{\roman{GL}}_{m_r})={\bold Z}^r\ \ .\leqno(4.7)$$ The element $\nu\in
X_M\otimes {\bold Q}$ in (ii) is of the form
$$\nu=(\nu(1)^{m_1},\ldots,
\nu(r)^{m_r})\ \ ,\ \hbox{with}\ m_i\cdot\nu(i)\in {\bold Z},\ \forall
i=1,\ldots, r\ \ .\leqno(4.8)$$ Then $\nu\in {\roman{Conv}}(W\mu )$ iff
$$0\leq \nu(i)\leq 1,\ \forall i=1,\ldots,r,\  \hbox{and}\
\Sigma_{i=1}^rm_i\cdot\nu(i)=s\ \ .\leqno(4.9)$$ It is therefore
obvious
that by letting $\nu$ vary over this set, its image in $\pi_1(M)$ is
equal
to $\kappa_M(W\mu )$. If $x=(s_1,\ldots, s_r)\in\kappa_M(W\mu
)$, then the
element $\tilde\nu$ is equal to $$\tilde\nu =(1^{s_1}, 0^{m_1-s_1},
1^{s_2}, 0^{m_2-s_2},\ldots, 1^{s_r}, 0^{m_r-s_r})\ \ ,$$ which
obviously
lies in $W\omega_s$.
\endbeweis

\bigskip\noindent {\bf Proposition 4.8.} {\it Let $G={\roman{GL}}_n$
or $G={\roman{GSp}}_{2n}$.
Let $\nu\in X_{M,{\bold Q}}\cap X_*(A)_{{\bold Q}, {\roman{dom}}}$
such
that its image under $\kappa_M\otimes {\bold Q}$ lies in
$\pi_1(M)$. Let
$\tilde\nu\in X_*(A)$ be the unique $M$-dominant $M$-minuscule
element
with $\kappa_M(\tilde\nu)= \kappa_M(\nu)$. Let $[\tilde\nu]$ be the
unique $G$-dominant
element in $W{\tilde\nu}$. Then we have, for every $\mu\in
X_*(A)_{{\roman{dom}}}$ with $\nu\leq\mu$, $$\nu\leq [\tilde\nu]\leq
\mu\
\ .$$ }

\medskip\noindent
{\bf Proof.} The first inequality is obvious since, $\nu$ being central
in
$M$, we have $$\nu\in{\roman{Conv}}(W_M\tilde\nu)\subset
{\roman{Conv}}(W\tilde\nu)\ \ .$$

Now we prove the second inequality. First we consider the case in which
$G={\roman{GL}}_n$, $M=M_{(m_1,\ldots,m_r)}$.
We note that
if $\mu$ is minuscule, then by the previous lemma we have
$\tilde\nu\in
W\mu$. Hence $[\tilde\nu]=\mu$, which proves the proposition in
this case.  We now
proceed by induction on $n$. As in the previous proof we write
$\nu=(\nu(1)^{m_1},\ldots, \nu(r)^{m_r})$.
\par
Assume first that there exists a maximal proper parabolic
$P'=M'N'$
containing $P$ such that $\nu\leq ^{M'}\mu$. This last condition
means
equivalently $$\nu\in {\roman{Conv}}(W_{M'}\mu)\Longleftrightarrow
\mu-\nu\in\sum_{\alpha^{\vee}\in\Delta_{M'}^{\vee}}{\bold
R}_+\alpha^{\vee}\ \ .$$ Here for any standard Levi subgroup $M$,
$\Delta_M^{\vee}$ denotes the set of simple coroots appearing in
$U\cap
M$.
\par
The Levi subgroup $M'$ corresponds to the partition
$(\sum_{j=1}^km_j,
\sum_{j=k+1}^r m_j)$ for some $k$ with $1\leq k\leq r-1$. Let us
subdivide
the interval $[0,n]$ into the $r$ subintervals $I(i)$, $i=1,\ldots, r$,
with $$I(1)= [0,m_1], I(1)=[m_1,m_1+m_2],\ldots, I(r)= [m_1+\ldots
+m_{r-1},n]\ \ .\leqno(4.10)$$ Since $\nu\leq^{M'}\mu$ we have
$\kappa_{M'}(\nu)=\kappa_{M'}(\mu)$ which means $$\sum_{j=1}^k
m_j\nu(j)=
\sum_{i\in I(1)\cup\ldots\cup I(k)}\mu_i\ \ ,\leqno(4.11)$$
where we adopt the convention that $\mu_0=0$ in order to make sense
of the right side of this
equation. The
converse
is also true by the following lemma applied to $M'=M$.
\medskip\noindent
{\bf Lemma 4.9.} {\it Let $\nu\leq\mu$ and
$\kappa_M(\nu)=\kappa_M(\mu)$.
Then $\nu\leq^M\mu$. }
\medskip\noindent
{\bf Proof.} We have by assumption
$$\mu-
\nu=\sum_{\alpha^{\vee}\in\Delta_G^{\vee}}c_{\alpha}\alpha^{\vee}\ \
,\ \ c_{\alpha}\in {\bold R}_+\ \ .$$ We want to show that
$c_{\alpha}=0$,
$\forall \alpha^{\vee}\in \Delta_G^{\vee}\setminus \Delta_M^{\vee}$.
But
$\Delta_G^{\vee}\setminus \Delta_M^{\vee}$ maps to a basis of
$X_{M,{\bold
Q}}/X_{G, {\bold Q}}$. Since $\kappa_M(\mu-\nu)=0$, we deduce
$$0=\kappa_M(\mu-\nu)\equiv
\sum_{\alpha^{\vee}\in\Delta_G^{\vee}\setminus
\Delta_M^{\vee}}c_{\alpha}\alpha^{\vee}{\roman{mod}}\ X_{G, {\bold
Q}}\ \
.\eqno\endbeweis$$ We will also need the following lemma.
\medskip\noindent
{\bf Lemma 4.10.} {\it Let $\mu\in X_*(A)_{\roman{dom}}$. Let
$\nu\in
X_*(A)$ be $M$-dominant with $\nu\leq^M\mu$. Then $\nu$ is
$G$-dominant and
  $\nu\leq \mu$. }
\medskip\noindent
{\bf Proof.} We have
$$\mu-
\nu=\sum_{\alpha^{\vee}\in\Delta_M^{\vee}}c_{\alpha}\alpha^{\vee}\ \
,\ \ c_{\alpha}\in {\bold R}_+\ \ .$$ Let $\beta\in \Delta_G\setminus
\Delta_M$. Then $\langle \alpha^{\vee},\beta\rangle \leq 0$, $\forall
\alpha^{\vee}\in\Delta_M^{\vee}$. Hence
$$\langle \nu,\beta\rangle
=\langle \mu,\beta\rangle
-\sum_{\alpha\in\Delta_M^{\vee}}c_{\alpha}\cdot\langle
\alpha^{\vee},\beta\rangle\geq 0\ \ .\eqno$$
Finally, it follows trivially from $\nu\leq^M\mu$ that $\nu\leq\mu$.
\endbeweis
\medskip\noindent
We apply this lemma to $M'$ and the unique $M'$-dominant
element
$[\tilde\nu]_{M'}$ in $W_{M'}\cdot\tilde\nu$. It satisfies
$[\tilde\nu]_{M'}\leq^{M'}\mu$ by induction hypothesis. By the
previous
lemma we therefore have $[\tilde\nu]_{M'}=[\tilde\nu]\leq \mu$ which
proves the proposition in this case.
\par
Let us now assume that there is no proper Levi subgroup $M'$
containing
$M$ such that $\nu\leq^{M'}\mu$. By Lemma 4.9 this means that
for
$k=1,\ldots,r-1$ we have $$\sum_{j=1}^km_j\cdot\nu(j) < \sum_{i\in
I(1)\cup\ldots\cup I(k)}\mu_i\ \ .\leqno(4.12)$$ In this case we are
going
to prove the assertion by induction on the height of $\mu$. If $\mu$
is
minuscule, the assertion is already proved. Otherwise there exists a
positive coroot $\alpha^{\vee}$ such that $\mu'= \mu-\alpha^{\vee}$
is
dominant. It suffices to show that $\nu\leq \mu'$ since then by
induction
hypothesis we have $[\tilde\nu]\leq\mu'$ and hence a fortiori
$[\tilde\nu]\leq\mu$.
\par
To prove $\nu\leq\mu'$ we introduce the partial sum functions for
$i=0,\ldots, n$, $$\leqalignno{ N_i &
=\sum_{\ell=1}^i\nu_{\ell}=\langle
\nu, \omega_i\rangle ,\ \  M_i= \sum_{\ell=1}^i\mu_{\ell}
=\langle\mu,
\omega_i\rangle\ \ , & (4.13)\cr M'_i & =\sum_{\ell=1}^i\mu'_{\ell}
=\langle \mu', \omega_i\rangle\ \ . \cr}$$ We have to show that
$N_i\leq
M'_i$, $\forall i=1,\ldots,n$, knowing that $N_i\leq M_i$, $\forall
i=1,\ldots, n$. These functions of $i$ may be interpolated into
continuous
functions on $[0,n]$ which are affine-linear on consecutive intervals
$[0,1], [1,2],\ldots$ and which are convex, since $\nu,\mu$ and
$\mu'$ are
all dominant. Furthermore, the function $N$ is affine-linear on the
intervals $I(1), I(2),\ldots, I(r)$. Hence it suffices to check that
$$N(x)\leq M'(x)\leqno(4.14)$$
for the endpoints $x$ of the intervals
$I(1),I(2),\ldots, I(r)$.
At the left endpoint of $I(1)$ and right endpoint of $I(r)$ we have
equality in (4.14). Now
consider the remaining endpoints. By (4.12) we have $N(x)<M(x)$. Since both
arguments are integers we conclude that $N(x)\leq M(x)-1$. On the
other
hand, since the positive coroot $\alpha^{\vee}$ is of the form
$\alpha^{\vee}=e_i-e_j$ for $i<j$ (where $e_1,\ldots, e_n$ is the
natural
basis of $X_*(A)={\bold Z}^n)$, it is obvious that $M'(x)\geq M(x)-1$,
which proves (4.14) in this case. This completes the proof for $\roman{GL}_n$.
  \par Now let us assume that
$G={\roman{GSp}}_{2n}$ and $M$ is the Levi subgroup obtained as the
intersection (inside $\roman{GL}_{2n}$) of ${\roman {GSp}}_{2n}$ and
$M_{(m_1,\ldots, m_r,2j,m_r,\ldots,m_1)}$. The second equality that
we need to prove for $G$ and
$M$ follows from that same inequality for $\roman{GL}_{2n}$ and
$M_{(m_1,\ldots,m_r,2j,m_r,\ldots,m_1)}$. To see this one simply
checks that ${\roman {GSp}}_{2n}$
inherits all relevant concepts (minuscule, dominant, $\leq$) from
$\roman{GL}_{2n}$, and
that the same is true for the two Levi subgroups).
\endbeweis
\medskip
Another way of formulating the previous proposition is that
$[\tilde\nu]$
is the unique minimal element of the set $$\{ \mu\in
X_*(A)_{\roman{dom}};\ \nu\leq \mu\}\ \ .\leqno(4.15)$$
\par
We can now prove the main result of this section.
\medskip\noindent
{\bf Theorem 4.11.} {\it Let $G={\roman{GL}}_n$ or
$G={\roman{GSp}}_{2n}$.
Let $b\in G(L)$, with associated Newton point
$\nu=\overline\nu(b)\in
X_*(A)_{{\bold Q},{\roman{dom}}}$. Let $\mu\in
X_*(A)_{\roman{dom}}$ with
$\nu\leq \mu$. Then the $\sigma$-conjugacy class of $b$ in $G(L)$
meets
$\tilde K\pi^{\mu}\tilde K$. Equivalently, there exists $h\in G(L)/\tilde
K$ with ${\roman{inv}}(h, b\sigma(h))=\mu$. }
\medskip\noindent
{\bf Proof.} After replacing $b$ by a $\sigma$-conjugate, we may
assume
that $b\in M(L)$ is basic, for a standard Levi subgroup $M$ [K I],
6.2. By
Proposition 4.6 we have to show that $\kappa_M(b)$ lies in the
image of
${\Cal P}_{\mu}$ in $\pi_1(M)$. By Proposition 4.8. we find
$\tilde\nu\in
{\Cal P}_{\mu}$ with
$\kappa_M(\tilde\nu)=\kappa_M(b)$.\endbeweis
\medskip
By Mazur's inequality we may summarize the previous theorem as
an equality
of two subsets of $X_*(A)_{\roman{dom}}$: Given $b\in G(L)$ we
have
$$\leqalignno{ & \{ \mu\in X_*(A)_{\roman{dom}};\ \exists h\in
G(L)/\tilde
K\ \hbox{with}\ {\roman{inv}}(h,b\sigma(h))=\mu\} & (4.16) \cr
=
& \{ \mu\in X_*(A)_{\roman{dom}};\ \overline\nu (b)\leq \mu\}\ \ .
\cr}$$
Furthermore, by (4.15) this subset has a unique minimal element.
\medskip\noindent
{\bf Remark 4.12.} Let $b\in M(L)$ be basic such that $M$ is the
centralizer of $\nu=\overline{\nu}(b)$ (i.e.\ $b$ is $G$-regular,
comp.\
[K I], 6.2.; recall that for any $\sigma$-conjugacy class $[b]$ there
exists a standard Levi subgroup $M$ and an element $b\in [b]$
which
satisfies these conditions). Let $\mu\in X_*(A)_{\roman{dom}}$
such that
$mb\sigma(m)^{-1}\in\tilde K_M \pi^{\mu} \tilde K_M$ for
some
$m\in M(L)$. It follows that $\nu\leq^M\mu$, or, equivalently by
Lemma
4.9, that $\nu\leq\mu$ and $\kappa_M(\nu)=\kappa_M(\mu)$.
Conversely let
$\mu\in X_*(A)_{\roman{dom}}$ such that $\nu\leq^M\mu$. Assume
furthermore
that there exists $g\in G(L)$ with $g^{-1}
b\sigma(g)\in\tilde
K\pi^{\mu}\tilde K$ and fix such an element. {\it Then, if
$G={\roman{GL}}_n$ or $G={\roman{GSp}}_{2n}$, it follows that
$g\in M(L)\cdot\tilde K$.}
\par
Indeed, assume that $G={\roman{GL}}_n$ and $M=M_{(m_1,\ldots,
m_r)}$. The
isocrystal $(N,F)=(L^n, b\sigma)$ has slope vector
$$\nu=(\nu(1)^{m_1},\ldots, \nu(r)^{m_r})\ \ ,$$ where
$\nu(1)>\nu(2)>\ldots >\nu(r)$. This chain of inequalities follows from
the assumption that $M$ is the centralizer of $\nu$ (equivalently,
the
break points of the Newton polygon of $(N,F)$ occur at
$m_1,m_1+m_2,\ldots, m_1+\ldots + m_{r-1})$. On the other hand,
$g\tilde
K$ defines a lattice $\Lambda$ in $N$ such that
$\mu(\Lambda)={\roman{inv}}(\Lambda,F\Lambda)=\mu$. Now the assumption
$\nu\leq^M\mu$ tells
us that the Hodge polygon of $\Lambda$ goes through all break points of
the
Newton polygon. Hence, by the Hodge-Newton decomposition
[Ka], Thm.\
1.6.1., we can write $$\Lambda=\bigoplus_{i=1}^r\Lambda_i\ \ ,\leqno(4.17)$$
where
$\Lambda_i=\Lambda\cap N_i$ is the intersection with the isotypic component
of slope
$\nu(i)$ of $N$. If the lattice $\Lambda_i$ corresponds to $g_i\cdot
{\roman{GL}}_{m_i}(O_L)$, then $g\tilde K=m\cdot \tilde K$
where
$m=(g_1,\ldots, g_r)\in\prod_{i=1}^r {\roman{GL}}_{m_i}(L)= M(L)$
which
proves the claim in this case. The case where
$G={\roman{GSp}}_{2n}$ is
similar.
\par
It seems likely that the above conclusion holds for more general
groups.
But we point out that the assumption that $\mu$ be $G$-dominant
is
essential; it is not enough to merely assume that $\mu$ is
$M$-dominant
with $\nu\leq^M\mu$, as the following example shows.
\par
Let $G={\roman{GL}}_3$, $M=M_{(2,1)}$ and $$b=\pmatrix
0&\pi^a&0\cr
\pi^{a+1}&0&0\cr 0&0&\pi^a\cr\endpmatrix\ \ ,\leqno(4.18)$$ where
$a>0$ is
a fixed integer. Then $\nu=\overline\nu(b)= (a+{1\over 2}, a+{1\over
2},
a)$. Let $\mu=(2a+1,0,a)$. Then $\mu$ is $M$-dominant but not
$G$-dominant
and $\nu\leq^M\mu$. We claim that there exists an element $g\in
G(L)\setminus M(L).\tilde K$ with $g^{-1}b\sigma(g)\in \tilde
K\pi^{\mu}\tilde K$. Let $$b'=\pmatrix0&\pi^a&0\cr \pi^{a+1}&0&0\cr
0&1&\pi^a\cr\endpmatrix\ \ .\leqno(4.19)$$ Then
$\overline\nu(b')=\nu$,
hence $b'$ is $\sigma$-conjugate to $b$ and ${\roman{inv}}(O_L^3,
b'\sigma(O_L^3))= [\mu]$, hence $b'\in\tilde K\pi^{\mu}\tilde
K$.
But an element $g\in G(L)$ with $g^{-1}b\sigma(g)=b'$ lies in
$M(L)\tilde
K$ if and only if $O_L^3$ is decomposable with respect to
the slope
decomposition of $L^3$ for $b'\sigma$. It is easy to see that $O_L^3$ is not decomposable.
\bigskip\noindent
{\bf 5. Restriction of scalars.} \medskip\noindent Let $F'$ be an
unramified field extension of degree $f$ of $F$. Let $V$ be a $F'$-
vector
space of dimension $n$. In the first part of this section we will be
concerned with the group $G=R_{F'/F}({\roman{GL}}(V))$. Let $\tilde
K\subset G(L)$ be a special maximal parahoric subgroup defined
over $F$.
The coset space $G(L)/\tilde K$ can be described as follows.
\par
We fix an embedding $F'\to L$. Then we can write
$$V\otimes_FL=\bigoplus\limits_{j\in {\bold Z}/f{\bold Z}}N_j\ \
,\leqno(5.1)$$ with $N_j=\{ v\in V\otimes_FL;\ (x\otimes 1)\cdot
v=(1\otimes \sigma^{-j}(x))\cdot x,\ \forall x\in F'\}$.
\par\noindent
Each summand is an $L$-vector space of dimension $n$. The
coset space
$G(L)/\tilde K$ parametrizes lattices for $O_{F'}\otimes_{O_F}O_L$
in
$V\otimes_FL$, or equivalently ${\bold Z}/f{\bold Z}$-graded
$O_L$-lattices, $$\tilde M=\bigoplus_{j\in {\bold Z}/f{\bold Z}}\tilde
M_j\ \ ,\leqno(5.2)$$ where each $\tilde M_j$ is an $O_L$-lattice in
$N_j$.
\par
Next we fix a conjugacy class of one-parameter subgroups of $G$.
Under the
decomposition $$G\otimes_FL= \prod_{j\in {\bold Z}/f{\bold
Z}}{\roman{GL}}(N_j)\leqno(5.3)$$ this corresponds to an $f$-tuple of
dominant cocharacters of ${\roman{GL}}_n$, $$\boldsymbol\mu
=(\mu_j)_{j\in {\bold Z}/f{\bold Z}}\ \ ,\ \ \mu_j\in ({\bold
Z}^n)_+\ \ ,\ \ \forall j=1,\ldots, f\ \ .\leqno(5.4)$$
\par
Finally, let $b\in G(L)$ and consider the $\sigma$-linear operator
$$F=b\cdot ({\roman{id}}_V\otimes \sigma)\leqno(5.5)$$ on
$V\otimes_FL$.
Then ${\roman{deg}}\ F=1$ with respect to the grading (5.1). We
introduce
the set $$\mathop{X}\limits^{\circ} (\boldsymbol\mu, b)_K=\{ \tilde
M=\bigoplus\limits_{j\in{\bold Z}/f{\bold Z}}\tilde M_j;\ \ {\roman{inv}}\
(\tilde M, F\tilde M) =\boldsymbol\mu\}\ \ .\leqno(5.6)$$ The last
condition is equivalent to ${\roman{inv}}(\tilde M_j, F\tilde
M_{j-1})=\mu_j$, $\forall j\in {\bold Z}/f{\bold Z}$, where each
invariant
is considered as an element of $({\bold Z}^n)_+$. Note that if
$\mu_1,\ldots,\mu_f$ are all minuscule, then by the minimality of
minuscule elements with respect to the partial order
$\mathop{\leq}\limits^!$ on $({\bold Z}^n)_+$, the set
$\mathop{X}\limits^{\circ}(\boldsymbol\mu, b)_K$ coincides with the set $X(\mu,
b)_K$
of section 3.
\medskip\noindent
{\bf Theorem 5.1.} {\it We have $$\mathop{X}\limits^{\circ}(\boldsymbol\mu,b)_K
\neq \emptyset\Longleftrightarrow [b]\in B(G,\boldsymbol\mu)\ \ .$$ \par}
\noindent
We note that the direct implication is
just
the group theoretic version of Mazur's theorem which was proved in
[RR].
To make the set $B(G,\boldsymbol\mu)$ more explicit, we note
the
Shapiro bijection [K II], 6.5.3. $$B(G,\boldsymbol\mu)= B(G',
\mu')\ \
.\leqno(5.7)$$ Here $G'={\roman{GL}}(V)$ is defined over $F'$ and
$$\mu'=\sum_j\mu_j\ \ .\leqno(5.8)$$ The map is obtained by
associating to
the $\sigma$-linear operator (5.5) of $V\otimes_FL$ the $\sigma^f$-
linear
operator on $N_0=V\otimes_{F'}L$, $$F^f:N_0\longrightarrow N_0\ \
.\leqno(5.9)$$ The condition that $[b]\in B(G,\boldsymbol\mu)$ is
equivalent to the condition that the slope vector $\boldsymbol\nu=
\boldsymbol\nu(F^f)\in ({\bold Q}^n)_+$ be smaller than
$\mu'$.
Let us now fix $b\in G(L)$ satisfying this condition and let us
construct
an element in $\mathop{X}\limits^{\circ}(\boldsymbol\mu, b)_K$.
Let
$\tilde M=\bigoplus_j\tilde M_j$ be any ${\bold Z}/f{\bold Z}$-graded
lattice and put $$M_j=F^j\tilde M_{-j}\ \ ,\ \ j=0,\ldots, f\ \
.\leqno(5.10)$$ Then $M_j$ is a lattice in $N_0$ and we obtain the
following description of $\mathop{X}\limits^{\circ}(\boldsymbol\mu,
b)_K$: $$\leqalignno{ \mathop{X}\limits^{\circ}(\boldsymbol\mu,b)_K= &
\{ (M_0, M_1,\ldots,M_f);\ M_f=F^fM_0, & (5.11) \cr &
{\roman{inv}}(M_j,
M_{j+1})=\mu_j,\ \forall j=0,\ldots, f-1\} \cr}$$
\medskip\noindent
We now apply Theorem 4.11 to the isocrystal $(V_0, F^f)$. Since
$\boldsymbol\nu\leq \mu'$ we obtain the existence of a lattice
$M_0$
in $N_0$ such that $${\roman{inv}}(M_0, F^fM_0)=\mu'\ \
.\leqno(5.12)$$ We
put $M_f=F^fM_0$. To complete the proof of Theorem 5.1, it
therefore
remains to fill in the remaining lattices $M_1,\ldots, M_{f-1}$. That
this
can be done follows from the following well-known lemma.
\medskip\noindent
{\bf Lemma 5.2.} {\it Let $\mu_0,\ldots, \mu_{f-1}\in ({\bold Z}^n)_+$
be
dominant vectors and let $\mu=\sum \mu_j$. Let $M_0, M_f$ be
lattices with
${\roman{inv}}(M_0, M_f)=\mu$. Then there exists a collection of
lattices
$M_1,M_2,\ldots, M_{f-1}$ such that ${\roman{inv}}(M_j,
M_{j+1})=\mu_j$,
$\forall j=0,\ldots, f-1$.\endbeweis }
\bigskip Now let $(V,\langle\ ,\
\rangle)$ be a symplectic vector space of dimension $2n$ over $F'$
and let
$G=R_{F'/F}({\roman{GSp}}(V,\langle\ ,\ \rangle))$. Let $\tilde
K\subset
G(L)$ be a special maximal parahoric subgroup defined over $F$.
The coset
space $G(L)/\tilde K$ parametrizes ${\bold Z}/f{\bold Z}$-lattices
$\tilde
M=\bigoplus_{j\in {\bold Z}/f{\bold Z}}\tilde M_j$ which are selfdual
with
respect to $\langle\ ,\ \rangle\otimes L$ up to a scalar in
$F'\otimes_FL$. Since the summands in (5.1) are orthogonal to
one another,
we may write $$\tilde M= \bigoplus_j\tilde M_j\ ,\ \hbox{where}\
\tilde
M_j^{\perp}=c_j\cdot \tilde M_j,\ c_j\in L,\ \forall j\in {\bold
Z}/f{\bold Z}\ \ .\leqno(5.13)$$
\par
A conjugacy class of one-parameter subgroups of $G$
corresponds to an
$f$-tuple of dominant cocharacters of ${\roman{GSp}}_{2n}$,
$$\boldsymbol\mu= (\mu_j)_{j\in {\bold Z}/f{\bold Z}}\ \ ,\ \ \mu_j\in
X_*(A)_{\roman{dom}}\ \ .\leqno(5.14)$$ Here $X_*(A)$ denotes the
cocharacter module for ${\roman{GSp}}_{2n}$ and
$X_*(A)_{\roman{dom}}=X_*(A)\cap ({\bold Z}^{2n})_+$.
\par
Finally, let $b\in G(L)$, with associated $\sigma$-linear operator
$F=b\cdot ({\roman{id}}_V\otimes \sigma)$ on $V\otimes_FL$. We
introduce
the set $$\leqalignno{ \mathop{X}\limits^{\circ}(\boldsymbol\mu,
b)_K=
\{ \tilde M=\bigoplus_{j\in{\bold Z}/f{\bold Z}}\tilde M_j;\ & \tilde
M_j^{\perp} =c_jM_j\ \hbox{for some}\ c_j\in L,\ \forall j, & (5.15) \cr
&
{\roman{inv}} (\tilde M_j, \tilde M_{j+1})=\mu_j,\ \forall j\}\ . \cr}$$
We introduce as before the lattices $M_j=F^j\tilde M_{-j}$ for
$j=0,\ldots, f$. Then, since $b$ is a symplectic similitude, it follows
that each lattice $M_j$ is self-dual up to a scalar. We therefore may
identify $\mathop{X}\limits^{\circ}(\boldsymbol\mu,b)_K$ with
  $$\leqalignno{
\quad
  \mathop{X}\limits^{\circ}(\boldsymbol\mu, b)_K= \{
  (M_0,\ldots,M_f);\ &
  M_f=F^fM_0,\ M_j^{\perp}=c_jM_j, \ \forall j,
  &
  (5.16)
  \cr
  &
  {\roman{inv}}(M_j, M_{j+1})= \mu_j,\ \forall j=0,\ldots, f-1\}\ \ .
  \cr}$$
  \medskip\noindent
  {\bf Theorem 5.3.} {\it We have
  $$\mathop{X}\limits^{\circ}(\boldsymbol\mu,b)_K\neq\emptyset\Longleftrightarrow [b]\in B(G,
  \boldsymbol\mu)\
  \ .$$
  }
  \medskip\noindent
  We only sketch the proof which is analogous to the proof of
Theorem 5.1.
  Again the direct implication follows from [RR]. To see the reverse
  implication, let us assume that $[b]\in B(G,\boldsymbol\mu)$, or
  equivalently, that the slope vector of $F^f:N_0\to N_0$ is smaller
  than $\mu'=\sum_j\mu_j$. An application of Theorem 4.11 shows
that there
  exists a lattice $M_0$ which is selfdual up to a scalar such that
  ${\roman{inv}}(M_0, F^fM_0)=\mu'$. We put $M_f=F^fM_0$.
Applying Lemma 5.2 we find
  a chain of lattices $M_0,M_1,\ldots,M_{f-1}$ such that
${\roman{inv}}(M_j,
  M_{j+1})=\mu_j$ for $j=0,\ldots,f-1$. But $M_0$ is selfdual up to a
  scalar and $\mu_0,\ldots, \mu_{f-1}\in X_*(A)_{\roman{dom}}$; this
implies
  successively that $M_1,M_2,\ldots,M_{f-1}$ are all selfdual up to a
  scalar. Hence we have indeed found an element of
  $\mathop{X}\limits^{\circ}(\mu, b)_K$.\endbeweis
  \bigskip\noindent
  {\bf 6. An incidence variety}
  \medskip\noindent
  Let $k$ be an algebraically closed field of characteristic $p$. We
fix a positive integer $f$.
  For each $i\in {\bold Z}/f{\bold Z}$ we fix a vector space $W_i$, all
of the
  same dimension $m>0$. Furthermore, for each $i\in{\bold Z}/f{\bold
Z}$ we fix
  a semi-linear map $\varphi_i:W_{i-1}\to W_i$ with respect to some
  automorphism $\sigma_i$ of $k$ and a semi-linear map
$\psi_i:W_i\to
  W_{i-1}$ with respect to some automorphism $\tau_i$ of $k$. We
assume that $\sigma_i$ and $\tau_i$
  are all powers (positive, negative, or zero) of the Frobenius
  automorphism of $k$.
  We impose the
  conditions
  $$\psi_i\circ\varphi_i=0\ \ ,\ \ \varphi_i\circ \psi_i=0\ \ ,\ \ \forall
  i\in {\bold Z}/f{\bold Z}\ \ .\leqno(6.1)$$
  We might picture these data in a circular diagram. Whenever you
turn back
  while traveling through this diagram you are killed (Orpheus
condition).
  $$\matrix
  &&
  W_1
  \cr\cr
  &
  \llap{$\scriptstyle\varphi_1$}\nearrow
\swarrow\rlap{$\scriptstyle\psi_1$}
  &&
\llap{$\scriptstyle\psi_2$}\nwarrow\searrow\rlap{$\scriptstyle
\varphi_2$}
  \cr\cr
  W_0&&&&W_2
  \cr\cr
  \llap{$\scriptstyle\varphi_0$}\big\uparrow\big\downarrow\rlap{$\scriptstyle\psi_0$}
  &&&&
  \big\uparrow\big\downarrow
  \cr\cr
  \vdots
  &&&&
  \vdots
  \cr\endmatrix
  \leqno(6.2)$$
  The aim of the present section is to prove the following theorem.
  \medskip\noindent
  {\bf Theorem 6.1.} {\it There exists a collection of lines
$\ell_i\subset
  W_i$ $(i\in {\bold Z}/f{\bold Z})$ such that
  $$\varphi_i(\ell_{i-1})\subset \ell_i\ \ ,\ \ \psi_i(\ell_i)\subset
  \ell_{i-1}\ \ ,\ \ \forall i\in {\bold Z}/f{\bold Z}\ \ .$$
  }
  \medskip\noindent
  Before starting the proof we make some comments. In the case
  $f=1$ we are given a vector space $W\neq (0)$ and two semi-linear
  endomorphisms $\varphi$ and $\psi$ of $W$ such that
  $\varphi\psi=\psi\varphi=0$. In this case we are looking for a line
  $\ell$ in $W$ which is carried into itself under $\varphi$ and $\psi$.
  This is essentially the situation considered in the proof of
Lemma 1.3
  where the existence of such a line is established. In the case
$f=2$ we
  are looking at a diagram
  $$
  \matrix
& \buildrel\varphi_1\over\longrightarrow \cr &
\buildrel\psi_1\over\longleftarrow \cr W_0 && W_1 \cr
&\buildrel\psi_0\over\longrightarrow \cr &
\buildrel\varphi_0\over\longleftarrow \cr\endmatrix
  \leqno(6.3)$$
  We are searching for a pair of lines $(\ell_0,\ell_1)$ which are
incident
  under $\varphi_0,\psi_0,\varphi_1,\psi_1$. In this case it is again
  possible to establish the existence of such a pair of lines by pure
  linear algebra. But already in this case a large number of case
  distinctions has to be made and this approach quickly gets out of
hand for a larger number of vector spaces.
Instead we use a density argument together with induction on $f$ to reduce the
problem to a special case that can be treated directly.
\medskip\noindent
{\bf Proof of Theorem 6.1.}
The special case goes as follows. Put
$$\Phi:=\varphi_f\varphi_{f-1}\dots\varphi_2\varphi_1,\leqno(6.4)$$
a semilinear endomorphism of~$W_0$, and assume that there
exists a line $\ell_0$ in $W_0$ such that
$\Phi \ell_0 =\ell_0$. For $i=1,\dots,f-1$ define a line $\ell_i$
in~$W_i$ by $\ell_i:=
\varphi_i\varphi_{i-1}\dots\varphi_2\varphi_1\ell_0$. Then
$\varphi_i\ell_{i-1}=\ell_{i}$ and
$\psi_i\ell_i=0$ for all $i \in \bold{Z}/f\bold{Z}$, so this
collection of lines solves our
problem.

The following reduction technique will be needed in the induction
on~$f$. Suppose $f>1$ and that
there exists $j$ such that $\psi_j$ is bijective, in which case
$\varphi_j$ is automatically $0$.
Given any family of lines $\ell_i$ solving our problem, we must have
$\ell_{j-1}=\psi_j\ell_j$.
Using $\psi_j$ to identify $W_{j-1}$ with $W_j$, and discarding the
two maps $\psi_j$,
$\varphi_j=0$, we are left with $f-1$ vector spaces
$\dots,W_{j-2},W_{j-1}\simeq
W_j,W_{j+1},\dots$ and maps $\varphi_i,\psi_i\quad(i\ne j)$ between them.
In other words we have a new problem of the same kind as our old one,
but with $f$ decreased
by~$1$. There is an obvious bijection between solutions of the old
and new problems.

The idea of the density argument is to fix $f$, $W_i$, $\sigma_i$,
$\tau_i$ and then to consider
the space~$M$ of all possible families of maps $\varphi_i,\psi_i$
satisfying condition $(6.1)$.
More precisely, for any finite dimensional $k$-vector spaces $W$,
$W'$ and any integral power
$\tau$ of Frobenius, we denote by $\operatorname{Hom}_\tau(W,W')$
the $k$-vector space of
$\tau$-linear maps $\varphi:W \to W'$, with scalar multiplication
by $\alpha \in k$ defined by
$(\alpha\varphi)(w)=\alpha(\varphi(w))$ (for all $w \in W$).
Returning to our fixed data
$f$, $W_i$, $\sigma_i$, $\tau_i$, we now put
$$
H_i:=\operatorname{Hom}_{\sigma_i}(W_{i-1},W_i) \times
\operatorname{Hom}_{\tau_i}(W_i,W_{i-1}),\leqno(6.5)
$$
a finite dimensional $k$-vector space which we regard as a
$k$-variety. Inside $H_i$ we have the
closed subvariety
$$
M_i:=\{(\varphi_i,\psi_i) \in H_i : \psi_i\varphi_i=0 \text{ and
$\varphi_i\psi_i=0$}\}.\leqno(6.6)
$$
Finally we put $M:=\prod_{i \in \bold Z/f\bold Z} M_i$, the space of
all families of maps
$\varphi_i$,$\psi_i$ satisfying $(6.1)$, which we are now regarding
as a (reducible) algebraic
variety over~$k$.

Writing $\bold P_i$ for the projective space of lines in~$W_i$, and
writing $\bold P$ for the
product $\bold P:=\prod_{i \in \bold Z/f\bold Z} \bold P_i$, we
consider the total incidence
variety $I \subset M \times \bold P$ consisting of
$(\varphi_i,\psi_i)_{i \in \bold Z/f\bold Z}
\in M$ and $(\ell_i)_{i \in \bold Z/f\bold Z} \in \bold P$ such that
$$
\varphi_i\ell_{i-1} \subset \ell_i \text{ and $\psi_i\ell_i \subset
\ell_{i-1}$}\leqno(6.7)
$$
for all $i \in \bold Z/f\bold Z$. It is easy to see that $I$ is
closed in $M \times \bold P$,
hence that the projection map $\pi:I \to M$ is proper. Thus
$M':=\pi(M)$ is closed in $M$, and
since Theorem 6.1 can be reformulated as the statement that $M'=M$,
it is enough to show that $M'$
is dense in~$M$. For this we need a better understanding of the
irreducible components of~$M$.

Recall that all the vector spaces $W_i$ have the same dimension
$m>0$. For any family $\bold
r=(r_i)_{i \in \bold Z/f\bold Z}$ of integers $r_i$ such that $0 \le
r_i \le m$, we denote by
$M_{\bold r}$ the subset of~$M$ consisting of families
$(\varphi_i,\psi_i)_{i \in \bold Z/f\bold Z}$
such that $\operatorname{rank}(\varphi_i)=r_i$ for all ${i \in \bold
Z/f\bold Z}$. (As usual
$\operatorname{rank}(\varphi_i)$ is the dimension of the image of
$\varphi_i$.) Thus $M$
has been decomposed into finitely many locally closed subsets
$M_\bold r$, and it is not hard to
see that each subset $M_\bold r$ is irreducible. (In the linear case,
{\it i.e.\/} when
$\sigma_i$, $\tau_i$ are the identity, the projection map
$(\varphi_i,\psi_i)_{i \in \bold Z/f\bold Z} \to (\varphi_i)_{i \in
\bold Z/f\bold Z}$
makes $M_\bold r$ into a vector bundle over a homogeneous space for a
product of general linear
groups; in general $M_\bold r$ is homeomorphic to such a vector
bundle, and is therefore
irreducible.)

For each $\bold r$ as above we are going to define a non-empty open
subset $U_\bold r$ of~$M_\bold
r$ such that $U_\bold r \subset M'$. This will show that $M'$ is
dense, as desired.

We define $U_\bold r$ to be the subset consisting of
$(\varphi_i,\psi_i)_{i \in \bold Z/f\bold Z}
\in M_\bold r$ satisfying the following two open conditions. The first is that
$\operatorname{rank}(\psi_i)=m-r_i$ for all $i$ (an open condition since
$\operatorname{rank}(\psi_i)\le m-r_i$ follows from $\varphi_i\psi_i=0$).

To formulate the second condition we again consider
$\Phi=\varphi_f\varphi_{f-1}\dots\varphi_2\varphi_1$, the semilinear
endomorphism of $W_0$ that appeared in our earlier discussion of
the special case of the theorem. Note that
$\operatorname{rank}(\Phi)
\le r_{\operatorname{min}}:=
\operatorname{min}\{r_i : i \in \bold Z/f\bold Z \}$. Let
$\Phi':\operatorname{im}\Phi \to
\operatorname{im}\Phi$ denote the restriction of $\Phi$ to the image
of $\Phi$ in~$W_0$. The
second condition is that
$\operatorname{rank}(\operatorname{\Phi})=r_{\operatorname{min}}$ and
that
the map $\Phi'$ be invertible. This is again an open condition on
$(\varphi_i,\psi_i)_{i \in \bold
Z/f\bold Z}
\in M_\bold r$.

We claim that $U_\bold r$ is non-empty. Choose a basis in each vector
space $W_i$, so that we can
represent the semilinear maps $\varphi_i$, $\psi_i$ by $m \times m$
matrices. Write $E_s$ for the
$m \times m$ diagonal matrix $1^s0^{m-s}$ and $F_s$ for the $m \times
m$ diagonal matrix
$0^s1^{m-s}$. Then put $\varphi_i=E_{r_i}\sigma_i$,
$\psi_i=F_{r_i}\tau_i$. Clearly
$(\varphi_i,\psi_i)_{i \in \bold Z/f\bold Z}$ lies in $U_\bold r$.

It remains to check that $U_\bold r \subset M'$. In other words for
$(\varphi_i,\psi_i)_{i \in \bold Z/f\bold Z} \in U_\bold r$ we must
show that there exists a
solution to the problem of finding lines $\ell_i \subset W_i$ such that
$$
\varphi_i\ell_{i-1} \subset \ell_i, \quad \psi_i \ell_i \subset \ell_{i-1}.
$$
There are two cases.

Suppose first that $r_{\operatorname{min}}=0$, so that there exists
$j \in \bold Z/f\bold Z$ such
that $r_j=0$. Thus $\varphi_j=0$ and it follows from the first open
condition that $\psi_j$ is
bijective. In this case we are done by induction on~$f$, as discussed earlier.

Now suppose that $r_{\operatorname{min}}>0$. By the second open
condition $\operatorname{im}\Phi
\ne 0$ and the restriction of $\Phi$ to $\operatorname{im}\Phi$ is
bijective. Therefore there
exists a line $\ell_0$ in $\operatorname{im}\Phi \subset W_0$ such
that $\Phi \ell_0 =\ell_0$. From
the special case treated directly at the beginning of this proof we
know that suitable lines
$\ell_i$ do exist, and thus the proof the theorem is complete.
\endbeweis
\medskip\noindent
{\bf Remark 6.2} (O.\ Gabber): The conclusion of Theorem 6.1 is
not true
without any hypotheses on the algebraically closed field $k$ and
the
automorphisms $\sigma_i$ and $\tau_i$ of $k$.  Indeed, in the
case when
all $\psi_i$ are zero, the theorem asserts the existence of an
eigenvector
of the semi-linear map $\Phi=\varphi_{f}\varphi_{f-1}\cdots \varphi_2\varphi_1$.
However,
such an eigenvector need not exist in general.
\bigskip\noindent
{\bf 7. General parahoric subgroups}
\medskip\noindent
In this section we will prove Conjecture 3.1 in the cases when
$G=R_{F'/F}({\roman{GL}}_n)$ or
$G=R_{F'/F}({\roman{GSp}}_{2n})$, where
$F'$ is an unramified extension of $F$.
\par
We start with the first group. We recall some notation from section
5. Let
$F'$ be an unramified extension of degree $f$ of $F$. Let $V$ be a
$F'$-vector space of dimension $n$. After fixing an embedding
$F'\to L$,
we have a decomposition (comp.\ (5.1)), $$V\otimes_FL=
\bigoplus\limits_{j\in {\bold Z}/f{\bold Z}}N_j\ \ .\leqno(7.1)$$ Let
$\bar I\subset{\bold Z}/n{\bold Z}$ be a non-empty subset. As in
section 1
we denote by $I$ the inverse image of $\bar I$ in ${\bold Z}$. A
${\bold
Z}/f{\bold Z}$-{\it graded periodic lattice chain of type} $\bar I$ is a
set of ${\bold Z}/f{\bold Z}$-graded $O_L$-lattices, one for each $i\in
I$, $$\tilde M^i=\bigoplus\limits_{j\in {\bold Z}/f{\bold Z}}\tilde M_j^i\
\ .\leqno(7.2)$$ Here $\tilde M_j^i=\tilde M^i\cap N_j$. We require
that,
for fixed $j$ the lattices $M_j^i$ form a periodic lattice chain of type
$\bar I$ in $N_j$, in the sense of (1.4). We denote by $X^G_{\bar
I}$ the
set of ${\bold Z}/f{\bold Z}$-graded periodic lattice chains of type
$\bar
I$. The conjugacy classes of parahoric subgroups of $G(L)$ defined
over
$F$ are in one-to-one correspondence with the non-empty subsets
$\bar I$
of ${\bold Z}/n{\bold Z}$. If $\tilde K$ is of type $\bar I$ we may
identify $G(L)/\tilde K$ with $X^G_{\bar I}$.
\par
We fix integers $r_j$ with $0\leq r_j\leq n$, $\forall j\in {\bold
Z}/f{\bold Z}$. We denote by $\boldsymbol\mu=(\mu_j)_{j\in{\bold
Z}/f{\bold Z}}$ the corresponding minuscule dominant coweight of
$G$, with
$\mu_j=\omega_{r_j}$.
\par
Let $b\in G(L)$. Then $b$ defines the $\sigma$-linear operator
$F=b\cdot({\roman{id}}_V\otimes\sigma)$ on $V\otimes_FL$,
which is of
degree 1 for the grading (7.1). Taking into account the identification
of
the $\mu$-admissible set with the $\mu$-permissible set for
${\roman{GL}}_n$ (compare the end of section 3) we may identify
the set
$X(\mu,b)_K$ of (3.7) in case $K$ is of type $\bar I$ with the
following
set $$\leqalignno{ X(\boldsymbol\mu,b)_{\bar I}=\{ & (\tilde
M_j^i)_{j,i}\in X^G_{\bar I};\ \tilde M_j^i\supset F\tilde
M^i_{j-1}\supset \pi\tilde M_j^i, & (7.3)\cr & {\roman{dim}}_{\bold
F}\tilde M_j^i/F\tilde M^i_{j-1}=r_j,\ \forall j\in{\bold Z}/f{\bold Z},\
\forall i\in I\}\ \ .\cr }$$ For each $i\in I$ and each $j$ with $0\leq
j\leq f$ let us set $M_j^i=F^j\tilde M^i_{-j}$. Let $N=N_0$. Then for
fixed $j$, the lattices $(M^i_j)_i$ form a periodic lattice chain of type
$\bar I$ in $N$. In section 1 we denoted the set of periodic lattice
chains of type $\bar I$ in $N$ by $X_{\bar I}$. Let us continue to do
so.
We therefore obtain from an element of $X^G_{\bar I}$ an $f$-tuple
of
elements of $X_{\bar I}$. We see that in this way we may identify
$X(\boldsymbol\mu,b)_K$ with the following set $$\leqalignno{
X(\boldsymbol\mu,b)_{\bar I}=\{ (M_j^i)_j\in & X_{\bar I}^f;\
M_f^i=F^fM_0^i, & (7.4) \cr & M_j^i\supset M^i_{j+1}\supset \pi
M_j^i,\
{\roman{dim}}_{\bold F}M_j^i/M^i_{j+1}=r_j\ , \cr & \forall i\in I,\
j=0,\ldots,f-1\ \}\ \ . \cr}$$
\medskip\noindent
{\bf Theorem 7.1.} {\it Conjecture 3.1 holds for
$G=R_{F'/F}({\roman{GL}}(V))$.}
\medskip\noindent
We proceed as in the proof of Proposition 1.1. If $\bar I$ consists
of a
single element, the statement (i) in Conjecture 3.1 follows from
Theorem
5.1. Again in statement (ii) it suffices to deal with the case when
$K$ is
an Iwahori subgroup, and this is then reduced to proving the
surjectivity
of the map $$X(\boldsymbol\mu, b)_{\bar I}\longrightarrow
X(\boldsymbol\mu, b)_{\bar J}\leqno(7.5)$$ when $\bar
J\subset\bar I$
differ by one element. It therefore suffices to show the following
analogue of Lemma 1.3.
\medskip\noindent
{\bf Lemma 7.2.} {\it Consider a commutative diagram of inclusions
of
lattices in $N$, $$\matrix M_0 & \supset & M_1 & \supset & \ldots
&
\supset & M_{f-1} & \supset & M_f &
=
& F^fM_0 \cr \cup && \cup &&&& \cup && \cup \cr M'_0 & \supset
& M'_1 &
\supset & \ldots & \supset & M'_{f-1}& \supset & M'_f & = &
F^fM'_0 \cr
\cup && \cup &&&& \cup && \cup \cr \pi M_0& \supset & \pi M_1
& \supset &
\ldots & \supset & \pi M_{f-1} & \supset & \pi M_f &
=
& \pi F^fM_0 & , \cr\endmatrix$$ where $M_j\supset
M_{j+1}\supset \pi M_j$
and $M'_j\supset M'_{j+1}\supset \pi M'_j$ with
$${\roman{dim}}_{\bold
F}M_j/M_{j+1}={\roman{dim}}_{\bold F}\ M'_j/M'_{j+1}=r_j\ \
\hbox{for}\
j=0,1,\ldots,f-1\ \ .$$ Assume that $M_j\neq M'_j$ for one $j$, or
equivalently, for all $j$. Then there exists a chain of lattices
$$L_0\supset L_1\supset \ldots \supset L_{f-1}\supset
L_f=F^fL_0$$ with
the following properties
\par\noindent
a) $L_j\supset L_{j+1}\supset \pi L_j$ with ${\roman{dim}}_{\bold
F}L_j/L_{j+1}=r_j$, for $j=0,\ldots,f-1$
\par\noindent
b) $M'_j\subset L_j\subset M_j$ with ${\roman{dim}}_{\bold
F}L_j/M'_j=1$,
for $j=0,\ldots, f-1$. }
\medskip\noindent
{\bf Proof.} Consider for $j=0,\ldots,f-1$ the ${\bold F}$-vector space
$$W_j=M_j/M'_j\ \ .\leqno(7.6)$$ These vector spaces have all the
same
dimension $\geq 1$. The inclusions $M_{j+1}\subset M_j$, resp.\
multiplication by $\pi$ induce linear maps $\psi$ resp.\ $\varphi$,
$$W_j
\matrix\buildrel\varphi\over\longrightarrow\cr
\mathop{\longleftarrow}\limits_{\psi}\endmatrix W_{j+1}\qquad
j=0,\ldots,
f-2\ \ .\leqno(7.7)$$ Similarly $F^f:M_0\to M_{f-1}$ and $\pi\cdot
(F^f)^{-1}:M_{f-1}\to M_0$ induce $\sigma^f$-linear resp.\
$\sigma^{-f}$-linear
maps $\psi$ resp.\ $\varphi$ $$W_0
\matrix\buildrel\varphi\over\longleftarrow\cr
\mathop{\longrightarrow}\limits_{\psi}\endmatrix W_{f-1}\ \
.\leqno(7.8)$$
It is obvious that we obtain in this way a diagram of the form (6.2)
which
satisfies all hypotheses of Theorem 6.1. We infer the existence of
lines
$(\ell_j\subset W_j)_j$ which are incident under the system of maps
$\varphi$ and $\psi$. Let $L_j\subset M_j$ be the inverse image of
$\ell_j$, for $j=0,\ldots, f-1$. Then we obtain a chain $L_0\supset
L_1\supset\ldots\supset L_{f-1}\supset L_f=F^fL_0$ which has the
required
properties.\endbeweis
\medskip\noindent
{\bf Variant 7.3.} Let $F'$ be an unramified extension of degree $f$
of
$F$ and let $D$ be a division algebra with center $F'$. Let $V$ be a
$D$-vector space of dimension $m$. Let
$G=R_{F'/F}({\roman{GL}}_D(V))$.
Then $G$ is an inner form of $R_{F'/F}({\roman{GL}}_n)$, where
$n=md$ with
$d^2={\roman{dim}}_{F'}D$. Since $G$ is not quasisplit, Mazur's
inequality
and its converse do not apply directly. Still we will show by
reduction to
the case of $R_{F'/F}({\roman{GL}}_n)$ that conjecture 3.1 holds
also in
this case. To simplify notations let us restrict ourselves to the case
$f=1$, i.e., $F'=F$.
\par
Let $O_D$ be the maximal order in $D$. Let $\tilde F$ be an
unramified
extension of $F$ of degree $d$ in $D$. Then we may write $O_D$
as $$O_D=
O_{\tilde F}[\Pi]/(\Pi\cdot a =a^{\sigma^s}\cdot\Pi\ \ \forall a\in
O_{\tilde F}\ \ ,\ \  \Pi^d=\pi)\ \ .\leqno(7.9)$$ Here $\Pi$ is a
uniformizer of $O_D$ and $s$ is inverse to the invariant of $D$ in
${\bold
Z}/d{\bold Z}$. Let us fix an embedding $\tilde F\to L$. Then we
obtain an
eigenspace decomposition analogous to (7.1),
$$V\otimes_FL=\bigoplus_{j\in{\bold Z}/d{\bold Z}}N_j\ ß
.\leqno(7.10)$$
Each $L$-vector space $N_j$ is of dimension $n$ and
${\roman{deg}}\ \Pi=s$
with respect to this grading. Put $N=N_0$. Let $\tilde K\subset
G(L)$ be a
parahoric subgroup maximal among those defined over $F$. Then
$G(L)/\tilde
K$ parametrizes the lattices $M$ in $V\otimes_FL$ which are
$O_D\otimes_{O_F}O_L$-invariant. Such a lattice is a free module
of rank
$m$ over $O_D\otimes_{O_F}O_L$. We associate to $M$ the
periodic lattice
chain in $N$, $${\Cal L}(M)=\ldots \supset M_0\supset \Pi M_{-
s}\supset
\Pi^2M_{-2s}\supset\ldots\supset \Pi^dM_0=\pi M_0\supset \ldots\ \
.\leqno(7.11)$$ Here $M_j=M\cap N_j$, so that
$M=\bigoplus_jM_j$. Then
${\Cal L}(M)$ is a periodic lattice chain of type $\bar I=\{ 0,m,\ldots,
(d-1)m\}\subset {\bold Z}/n{\bold Z}$. In this way we obtain a
bijection
$$G(L)/\tilde K=X_{\bar I}\simeq {\roman{GL}}_n(L)/\tilde K_{\bar I}\ \
.\leqno(7.12)$$ Here $K_{\bar I}$ is a parahoric subgroup of
${\roman{GL}}_n$ defined over $F$ and we have implicitly chosen a
basis of
the $L$-vector space $N$. More generally, we obtain a bijection
between
the sets of conjugacy classes of parahoric subgroups of $G(F)$
and the
non-empty subsets $\bar I\subset {\bold Z}/n{\bold Z}$ which are
invariant
under the translation action $x\mapsto x+\overline m$ on ${\bold
Z}/n{\bold Z}$. If $K$ corresponds to $\bar I$, then again
$$G(L)/\tilde
K=X_{\bar I}\simeq {\roman{GL}}_n(L)/\tilde K_{\bar I}\ \ .$$
\par
Let $b\in G(L)$. Then $b$ defines the $\sigma$-linear operator
$F=b\cdot
({\roman{id}}_V\otimes \sigma)$ on $V\otimes_FL$. The relation
between the
Newton point of $b$ and the slope vector of $F$ is given by
  $$\boldsymbol\nu (F)=
(\nu(b)^d)\ \ .\leqno(7.13)$$ Here $\nu(b)^d\in ({\bold Q}^{nd})_+$ is
the
vector which repeats $d$ times each entry of $\nu(b)\in ({\bold
Q}^n)_+$.
Let $\tilde K\subset G(L)$ be a parahoric subgroup maximal among
those
defined over $F$. We identify $\tilde W^K\setminus \tilde W/\tilde
W^K$
with $\tilde W^{\bar I}\setminus \tilde W/\tilde W^{\bar I}$, where
$\bar
I=m\cdot {\bold Z}/n{\bold Z}$ and where $\tilde W^{\bar I}=\tilde
W^{\tilde K_{\bar I}}$. (Something analogous holds for any parahoric
subgroup of $G(L)$ defined over $F$). Let $g,g'\in G(L)/\tilde K$.
Let $M$
and $M'$ be the corresponding $O_D\otimes_{O_F}O_L$-stable
lattices in
$V\otimes_FL$, with corresponding decompositions $M=
\bigoplus_{j\in{\bold
Z}/d{\bold Z}}M_j$ and $M'=\bigoplus_{j\in {\bold Z}/d{\bold Z}}M'_j$
and
corresponding periodic lattice chains ${\Cal L}(M)$ and ${\Cal
L}(M')$.
Let $\mu\in ({\bold Z}^n)_+$ be a dominant cocharacter of $G$.
Then
$$\leqalignno{ {\roman{inv}}(g,g')\in {\roman{Adm}}_{\tilde K}(\mu) &
\Longleftrightarrow {\roman{inv}}({\Cal L}(M), {\Cal L}(M'))\in
{\roman{Adm}}_{\tilde K_{\bar I}}(\mu) & (7.14) \cr &
\Longleftrightarrow
{\roman{inv}} (M_j, M'_j)\leq\mu\ \ ,\ \ \forall j\in {\bold Z}/d{\bold
Z}\ \ . \cr}$$ Assume now that $g'=b\sigma(g)$. Comparing the
lattices $M$
and $M'$ in $V\otimes_FL$, we obtain from (7.14) that
$${\roman{inv}}(M,FM)\leq (\mu^d)\ \ .\leqno(7.15)$$ (Here
${\roman{inv}}$
denotes the relative position of two lattices in $V\otimes_FL\simeq
L^{nd}$.) By Mazur's inequality we infer that $\boldsymbol\nu(F)\leq (\mu^d)$.
By
(7.13) this implies that $\nu(b)\leq \mu$.
\par
This proves one implication of statement (i) in Conjecture 3.1. The
remaining assertions of the Conjecture follow from the case of
${\roman{GL}}_n$ and the preceding remarks connecting the case
at hand to
the case of ${\roman{GL}}_n$.\endbeweis
\medskip
We now turn to the case $G=R_{F'/F}({\roman{GSp}}_{2n})$. We
now let $V$
be an $F'$-vector space of dimension $2n$ equipped with a
non-degenerate symplectic form
$\langle\ ,\ \rangle$. Let $G= R_{F'/F}({\roman{GSp}}(V,\langle\ ,\
\rangle))$. The decomposition (7.1) is an orthogonal sum
decomposition
with respect to $\langle\ ,\ \rangle$. Hence each summand $N_j$
is a
symplectic vector space of dimension $2n$ over $L$.
\par
The conjugacy classes of parahoric subgroups of $G(L)$ defined
over $F$
are in one-to-one correspondence with the non-empty symmetric
subsets
$\bar I$ of ${\bold Z}/2n{\bold Z}$. A ${\bold Z}/f{\bold Z}$-graded
periodic lattice chain $(\tilde M^i_j)_{i,j}$ of type $\bar I$ is called
selfdual, if for each $j\in{\bold Z}/f{\bold Z}$ the periodic lattice
chain $\tilde M_j$ of type $\bar I$ is selfdual in the sense of (2.3).
We
denote by $X^G_{\bar I}$ the set of ${\bold Z}/f{\bold Z}$-graded
selfdual
periodic lattice chains of type $\bar I$. If $\tilde K$ is a parahoric
subgroup of type $\bar I$ defined over $F$, we may identify the
coset
space $G(L)/\tilde K$ with $X^G_{\bar I}$.
\par
We fix integers $r_j\in \{ 0,n,2n\}$, $\forall j\in{\bold Z}/f{\bold Z}$.
We denote by $\boldsymbol\mu=(\mu_j)_{j\in{\bold Z}/f{\bold Z}}$
the
corresponding minuscule dominant coweight of $G$, with
$\mu_j=\omega_{r_j}$.
\par
Let $b\in G(L)$. Then $b$ defines the $\sigma$-linear operator
$F=b\cdot
({\roman{id}}_V\otimes \sigma)$ on $V\otimes_FL$. It is of degree
1 with
respect to the grading (7.1) and there are scalars $c_j\in L$, $\forall
j\in {\bold Z}/f{\bold Z}$ such that $$\langle Fv, Fw\rangle
=c_j\cdot\langle v,w\rangle\ \ ,\ \ v,w\in N_j\ \ .\leqno(7.16)$$ We
have
the following description of the set $X(\mu,b)_K$ of (3.7) in case
$K$ is
of type $\bar I$, $$X^G(\boldsymbol\mu,b)_{\bar I}=\{ (\tilde
M_j^i)_{i,j}\in X^G_{\bar I};\ (\tilde M_j^i)\in X(\boldsymbol\mu,
b)_{\bar I}\}\ \ .\leqno(7.17)$$ Here $X(\boldsymbol\mu,b)_{\bar I}$
is the set (7.3) for the group $R_{F'/F}({\roman{GL}}(V))$. In other
words, the elements of $X^G(\boldsymbol\mu,b)_{\bar I}$ are the
${\bold Z}/f{\bold Z}$-graded selfdual periodic lattice chains
$(\tilde M_j^i)$ of type
$\bar I$ which satisfy $$\tilde M_j^i\supset F\tilde
M^i_{j-1}\supset \pi \tilde M^i_j\ \
;\ \ {\roman{dim}}\ \tilde M^i_j/F\tilde M^i_{j-1}=r_j\ ,\ \forall i\in I\ ,\
\forall j\in {\bold Z}/f{\bold Z}\ \ .\leqno(7.18)$$
\medskip\noindent
{\bf Theorem 7.4.} {\it Conjecture 3.1 holds for
$G=R_{F'/F}({\roman{GSp}}(V,\langle\ ,\ \rangle))$. }
\medskip
We proceed as in the proof of Proposition 2.1. If $\bar I=\{ 0\}$, the
statement (i) of Conjecture 3.1 follows from Theorem 5.3. In
statement
(ii) it suffices to deal with the case when $K$ is an Iwahori
subgroup,
and this is then reduced to proving the surjectivity of the map
$$X^G(\boldsymbol\mu,b)_{\bar I}\longrightarrow X^G(\boldsymbol\mu,b)_{\bar J}\ \ ,\leqno(7.19)$$
in the situation considered in
(2.5). In other words, $$\bar I= \bar J\cup \{ \overline k+\overline 1,
-(\overline k+\overline 1)\}\ \ ,\ \hbox{where}\ \overline k\in\overline
J\ \hbox{with}\ \overline k+\overline 1\not\in \bar J\ .$$ We choose
as in
(2.5) a representative $k$ of $\overline k$ in ${\bold Z}$ and denote
by
$\ell$ the next largest element in $J$. The assertion reduces to the
corresponding statement for ${\roman{GL}}_{2n}$ (Lemma 7.2) with
the
following variant of Lemma 2.3.
\medskip\noindent
{\bf Lemma 7.4.} {\it Let $\bar J\subset\bar I$ as above. Let $(\tilde
M_j^i)$ be a ${\bold Z}/f{\bold Z}$-graded selfdual periodic lattice
chain
of type $\bar J$. The set of refinements of $(\tilde M_j^i)$ into a
${\bold Z}/f{\bold Z}$-graded selfdual periodic lattice chain of type
$\bar I$ is in one-to-one correspondence with the set of ${\bold
Z}/f{\bold Z}$-graded lattices $(\tilde M_j)_{j\in{\bold Z}/f{\bold Z}}$
with the property that $$\tilde M_j^k\subset \tilde M_j\subset \tilde
M_j^{\ell}\ \ \hbox{and}\ \ {\roman{dim}}_{\bold F} \tilde M_j/\tilde
M_j^k=1,\ \ \forall j\in{\bold Z}/f{\bold Z}\ \ .$$ }
\medskip\noindent
{\bf Proof.} This follows from Lemma 2.3 applied to each direct
summand
$N_j$, $j\in{\bold Z}/f{\bold Z}$.
\bigskip\bigskip
\centerline{\bf References}
\bigskip\noindent
  \ref{[HN]} Haines, T., Ngo, B.C.: Alcoves associated to special
fibers of
  local models, preprint Toronto 2000.
  \ref{[Ka]} Katz, N.: Slope filtration of $F$-crystals. Ast\'erisque
  {\bf 63} (1979), 113--164.
   \ref{[K I]} Kottwitz, R.: Isocrystals
with additional structure. Compositio Math. {\bf 56} (1985), no. 2,
201--220. \ref{[K II]} Kottwitz, R.: Isocrystals with additional
structure. II. Compositio Math. {\bf 109} (1997), no. 3, 255--339.
\ref{[KR]} Kottwitz, R., Rapoport, M.: Minuscule alcoves for
${\roman{{\roman{GL}}}}\sb n$ and $G{\roman{Sp}}\sb {2n}$.
Manuscripta
Math. {\bf 102} (2000), no. 4, 403--428.
\ref{[LR]} Langlands, R., Rapoport, M.: Shimuravariet\"aten und Gerben. J.
Reine Angew. Math. {\bf 378} (1987), 113--220.
  \ref{[RR]} Rapoport, M.,
Richartz, M.: On the classification and specialization of $F$-
isocrystals
with additional structure. Compositio Math. {\bf 103} (1996), no. 2,
153--181.
  \ref{[R]} Rapoport,
M.: A positivity property of the Satake isomorphism. Manuscripta
Math.
{\bf 101} (2000), no. 2, 153--166.

\end